\documentclass[a4paper,12pt]{article}
\usepackage{graphicx}
\usepackage{amsmath,amssymb}
\usepackage{amsfonts}

\usepackage{color}

\newtheorem{thm}{Theorem}

\newtheorem{case}{Case}
\newtheorem{subcase}{Subcase}
\numberwithin{subcase}{case}

\newtheorem{la}[thm]{Lemma}

\newtheorem{cor}[thm]{Corollary}
\newtheorem{prop}[thm]{Proposition}

\newtheorem{prob}[thm]{Problem}

\newtheorem{alg}{Algorithm}[section]

\newcommand{\white}{\qquad{\framebox{\rule{0pt}{4pt}}}}
\newcommand{\owari}{\hfill\white}

\DeclareMathOperator{\Od}{O}

\DeclareMathOperator{\Eod}{\Theta}
\DeclareMathOperator{\IA}{IA}
\DeclareMathOperator{\IC}{IC}
\DeclareMathOperator{\Irr}{Irr}
\DeclareMathOperator{\id}{id}

\newcommand{\qedgo}{\vspace{0.35cm}}

\newcommand{\tume}{\hspace*{0.2em}}

 \def\np{\mathcal{NP}}

\def\co{\textrm{co-}}

\newcommand{\ol}[1]{\overline{#1}}

\newcommand{\maru}[1]{{\ooalign{\hfill$\scriptsize#1$\hfill\crcr$\bigcirc$}}}

\newenvironment{alg-enumerate}{%
  \begin{enumerate}%
 }{\end{enumerate}%
}

\newenvironment{u-enumerate}{%
  \begin{enumerate}%
 }{\end{enumerate}%
}

\newenvironment{Proof}{\expandafter
\paragraph{{\bf Proof\/.}}}{\hspace{5mm}\owari\qedgo}

\newenvironment{Proofof}[1]{\expandafter
\paragraph{{\bf Proof of #1.}}}{\hspace{5mm}\owari\qedgo}

  \newenvironment{namelist}[1]{%
    \begin{list}{}
      {
      \settowidth{\labelwidth}{#1}
      \setlength{\leftmargin}{1.1\labelwidth}
      \setlength{\labelsep}{5pt}
      \setlength{\parsep}{0pt}}
  }{%
  \end{list}}

\begin{document}
\begin{center}
{\Large Agent Arrangement Problem}
\end{center}
\begin{center}
{\large Tomoki Nakamigawa\footnote{Email: {\texttt{nakami@info.shonan-it.ac.jp}}}
}
\end{center}
\begin{center}
{
Department of Information Science \\%
Shonan Institute of Technology \\%
1-1-25 Tsujido-Nishikaigan, Fujisawa 251-8511, Japan
}
\end{center}
\begin{center}
{\large Tadashi Sakuma\footnote{This work was supported by Grant-in-Aid for Scientific Research (C).\\
Email: {\texttt{sakuma@e.yamagata-u.ac.jp}} 
}
}
\end{center}
\begin{center}
{
Systems Science and Information Studies \\%
Faculty of Education, Art and Science \\%
Yamagata University \\%
1-4-12 Kojirakawa, Yamagata 990-8560, Japan
}
\end{center}
\begin{abstract}
An {\em arrangement} of an ordered pair $(G_A, G_M)$ of graphs 
is defined as a function $f$ from $V(G_A)$ to $V(G_M)$ such that, 
for each vertex $c$ of  $G_M$, the vertex-set $f^{-1}(c)$ of $G_A$ 
either is $\emptyset$ (the case when $c \not\in f(V(G_A))$) or 
induces a connected subgraph of $G_A$ and that the family 
$\{f^{-1}(y) : y \in V(G_M), f^{-1}(y) \neq \emptyset\}$ 
is a partition of $V(G_A)$. Let $f$ be an arrangement of 
$(G_A, G_M)$, let $pq$ be an edge of $G_M$ and let $U$ be 
a subset of $f^{-1}(p)$ such that each of the three graphs 
$G_A[U]$, $G_A[f^{-1}(p)\setminus U]$ and $G_A[f^{-1}(q)\cup U]$ 
is either connected or $\emptyset$ and that 
$\big( f^{-1}(p)\cup f^{-1}(q) \big) \setminus U \neq \emptyset$. 
A {\em  transfer} of $U$ from $p$ to $q$ is defined as the 
modification $f^{\prime}$ of $f$ such that $f^{\prime}(x):=f(x)$ 
for every $ x \notin U$ and $f^{\prime}(u):=q$ for every $u \in U$. 
Two arrangements $f$ and $g$ of $(G_A, G_M)$ are called 
{\em  t-equivalent} if they can be transformed into each other 
by a finite sequence of transfers. An ordered pair $(G_A, G_M)$ of 
graphs is called {\em  almighty} if every two arrangements of 
the pair $(G_A, G_M)$ are t-equivalent. In this study, we consider 
the following two decision problems. 
\begin{itemize}
\item[{\bf (P1)}]{For a given pair of arrangements $f$ and $g$ of a 
given ordered pair $(G_A,G_M)$ of graphs, decide whether $f$ is
t-equivalent to $g$ or not. }
\item[{\bf (P2)}]{For a given ordered pair $(G_A,G_M)$ of graphs, 
decide whether the pair $(G_A,G_M)$ is almighty or not. }
\end{itemize}
We show an $\Od(|E(G_A)|+(|V(G_M)|+|E(G_A)|)|V(G_A)|)$-time algorithm for 
{\bf (P1)}. By using this algorithm, we can also construct an 
explicit sequence of transfers from $f$ to $g$ of
$\Eod(|V(G_M)|^2 \cdot |V(G_A)|)$-length. Lastly we prove 
the $\co\np$-completeness of {\bf (P2)}. 
\end{abstract}
keywords: pebble motion, motion planning

\section{Introduction}

In this paper, we introduce two generalizations of the pebble 
motion problems \cite{CDP2006, J1879, KMS1984, 
PRST1994, W1974}, the first of which we call ``the Subgraph 
Allocation Problem''. The second is a specialization of the first, 
which we call ``the Agent Arrangement Problem''.  We take up 
this specialization because it has some significance as a 
theoretical model of logistics.

Throughout this paper, a graph is undirected with no loop 
or multiple edge. For a graph $G$, let $V(G)$ and $E(G)$ 
denote the vertex set of $G$ and the edge set of $G$. 
Let $G_A$ be a simple undirected graph such that its 
vertex set $V(G_A)$ is the set of {\em agents} and that 
$G_A$ has an edge $uv$ if and only if the two agents 
$u$ and $v$ have a method of taking mutual 
communication. Let us call the graph $G_A$ an 
{\em agent  network}. On the other hand, let $G_M$ 
be a simple undirected graph such that its vertex set 
$V(G_M)$ is the set of {\em countries} and that $G_M$ 
has an edge $pq$ if and only if there exits a way of 
mutually direct transportation between the two 
countries $p$ and $q$. Let us call the graph $G_M$ 
a {\em route map}. An {\em arrangement} of the 
ordered pair $(G_A, G_M)$ of graphs is defined as 
a function $f$ from $V(G_A)$ to $V(G_M)$ such that 
the family of sets $\{f^{-1}(c) : c \in V(G_M), f^{-1}(c) \neq 
\emptyset\}$ is a partition of $V(G_A)$ and each set 
$f^{-1}(c)$ in the family (i.e. the agents staying in the country 
$c$) induces a connected subgraph of $G_A$. Let $f$ be an 
arrangement of $(G_A, G_M)$, let $st$ be an edge of $G_M$ 
and let $U (\subseteq f^{-1}(t))$ be a set of agents staying in 
the country $s$ such that each of the three graphs $G_A[U]$, 
$G_A[f^{-1}(s)\setminus U]$ and $G_A[f^{-1}(t)\cup U]$ is ether 
connected or $\emptyset$. 
A {\em  transfer} of the set of agents $U$ along with the edge 
$st$ of $G_M$ from the country $s$ to the country $t$ is 
defined as the modification $f^{\prime}$ of $f$ such that 
$f^{\prime}(x):=f(x)$ for every $ x \notin U$ and $f^{\prime}(u):=t$ 
for every $u \in U$. Note that the modification $f^{\prime}$ is 
again an arrangement of the pair $(G_A, G_M)$. 
Here we use the {\em Subgraph Allocation Problem} ({\em SGA}, 
for short)  as a general term for the set of related problems 
dealing with arrangements and agent-transfers on given ordered 
pairs of graphs. 

In addition to the previous definition of a 
{\em transfer} of a set of agents $U$ from a country $s$ 
to a country $t$, we sometimes need to assume that such 
a transfer is allowed if and only if at least one of the following 
two conditions hold true:
\begin{enumerate}
\item{$f^{-1}(t) \neq \emptyset$. }
\item{$f^{-1}(s) \setminus U \neq \emptyset$.}
\end{enumerate}
The meaning of these two additional conditions can be interpreted
as follows. 
The first condition describes the situation when the agent subnetwork 
$G_A[U]$ will move from the country $s$ to the country $t$ relying on 
the connection to the agent subnetwork $G_A[f^{-1}(t)]$ in the country 
$t$.  If $f^{-1}(t) = \emptyset$ and the `vanguard' agent subnetwork 
$G_A[U]$ of the agent organization $G_A$ will pioneer the new territory 
$t$, then $G_A[U]$ must need backup support of the remainder agent 
subnetwork $G_A[f^{-1}(s) \setminus U]$ in the home ground $s$. Hence 
the second condition is necessary if the first condition does not hold. 
Note that the above pair of conditions can be summarized as the single 
condition $f^{-1}(\{s, t\}) \setminus U \neq \emptyset$.
When we impose this restriction on a transfer, to avoid 
the confusion to the unrestricted version, we use the term 
the {\em Agent Arrangement Problem} ({\em AAP}, for short) 
instead of SGA. 
The Agent Arrangement Problem can be regarded as a natural 
model for logistics, that is, an appropriate treatment for 
(re-)configurations on (distribution and/or human) networks in 
security. And hence it has clear applications to the wide area such 
as computer science, engineering, social and political science. 

Note that both of the SGA and AAP models are generalizations of the 
pebble motion problem. Actually, on the SGA model, if we assume the 
agent network $G_A$ to be an edge-less graph, the set of agents 
$V(G_A)$ can be regarded as the set of pebbles on the board graph 
$G_M$. In this case, a transfer of an agent turns to be a move of a 
pebble on the board. In the same way, on the AAP model, if we assume 
the agent network $G_A$ to be a disjoint union of two-vertex complete 
graphs, and if we exclude from the consideration the meaningless 
arrangements $f$ which contain a `frozen' pair $\{u,v\}$ of agents 
(that is, an edge $\{u, v\}$ of $G_A$ such that $\{f(u), f(v)\} \notin 
E(G_M)$), then the behavior of the model is essentially the same as 
the behavior of the pebble motion problem. Of course, for general 
cases, the behavior of transfers on SGA or AAP is far from 
the behavior of pebble motions, and its analysis is considerably 
difficult issue as is shown in Section\tume\ref{almighty}. 

On both of the SGA and AAP models, two arrangements $f$ and $g$ 
of an ordered pair $(G_A, G_M)$ of graphs is called {\em  t-equivalent} 
and is denoted by $f \cong g$ if the two arrangements can be 
transformed into each other by a finite sequence of transfers. 
An ordered pair $(G_A,G_M)$ of graphs is called {\em  almighty} if 
every two arrangements of the pair $(G_A, G_M)$ are t-equivalent. 
In this study, we consider the following two decision problems. 
\begin{itemize}
\item[{\bf (P1)}]{For a given two arrangements $f$ and $g$ of a 
given ordered pair $(G_A,G_M)$ of graphs, decide whether $f$ is 
t-equivalent to $g$ or not. }
\item[{\bf (P2)}]{For a given ordered pair $(G_A,G_M)$ of graphs, 
decide whether the pair $(G_A,G_M)$ is almighty or not. }
\end{itemize}
First, by using the previous results of Wilson\tume\cite{W1974} and 
Kornhauser et al.\tume\cite{KMS1984},  we will give a proper description 
of several good characterizations to the equivalence decision 
(and some related problems) for the classical pebble motion 
problem, which we need in our polynomial algorithm for {\bf(P1)}. 
We provide these in Sections $3$ and $4$. 
In Section\tume\ref{equivalency}, we describe 
an$\Od(|E(G_M)|+(|V(G_M)|+|E(G_A)|)|V(G_A)|)$-time algorithm 
for {\bf (P1)} on AAP and SGA, in common. 
By using this algorithm, we can also construct an 
explicit sequence of transfers from $f$ to $g$ of
$\Eod(|V(G_M)|^2 \cdot |V(G_A)|)$-length.   
In Section\tume\ref{almighty}, we prove the $\co\np$-completeness of 
{\bf (P2)} for the both of SGA and AAP.


\section{Preliminary for the Pebble Motion Problem}
Let $G$ be a connected graph with $n$ vertices.
Let $P$ be the set of labeled pebbles of order $m < n$.
A \textit{configuration} of the set of pebbles $P$ on $G$ is 
defined as an injective function $f$ from $P$ to $V(G)$, 
where if $f^{-1}(v) \neq \emptyset$, then $f^{-1}(v)$ 
represents a pebble of $P$ on the vertex $v$ of $G$,  
and if $f^{-1}(v)=\emptyset$, it means that $v$ is unoccupied.
Any pebble $p \in P$ must be on some vertex of the graph.
Hence we have $|f^{-1}(v)| \leqq 1$ for any $v \in V(G)$.

A \textit{move} is transferring a pebble to an adjacent 
unoccupied vertex. For a pair of configurations $f$ and $g$, 
we say that $f$ and $g$ are \textit{equivalent} if $f$ can be 
transformed into $g$ by a sequence of finite moves. 
we write $f \sim g$ if $f$ and $g$ are equivalent. 

Let us define the puzzle graph ${\rm puz}(G,k)$ of a graph $G$ 
with $k$ unoccupied vertices such that $V({\rm puz}(G,k))$ is 
the set of all the configurations ${\cal F}(G)$, and 
$E({\rm puz}(G,k)) = \{ (f,g) \,:\, f,g \in {\cal F}(X), f$ can be 
transformed into $g$ by a single move$\}$. 
For example, if $G$ is a $4 \times 4$ grid graph, 
${\rm puz}(G,1)$ corresponds to a well-known 
``15 puzzle''\cite{A1999,J1879,S1879}.

We say that $(G,k)$ is \textit{transitive} if for any configuration 
$f$ and for any pebble $p$ of $P$, $p$ can be shifted to an 
arbitrary vertex of $G$ by a sequence of finite moves. For a graph 
$G$, let $c(G)$ be the number of connected components of $G$. 
We say that $(G,k)$ is \textit{feasible} if $c({\rm puz}(G,k)) = 1$.

Wilson studied the problem for the case $k=1$ \cite{W1974}. 
It is not difficult to see that $(G,1)$ is transitive if and only if 
$G$ is $2$-connected. Let $S_m$ and $A_m$ denote the 
symmetric group and the alternating group of order $m$, 
respectively. For a finite set $M$, let $S(M)$ be the symmetric 
group on $M$. For a vertex $x$ of $V(G)$, let ${\cal F}_x$ 
be the set of configurations $f$ with $f^{-1}(x)=\emptyset$ and 
define $G_x$ as the set of permutations $\sigma \in S(V(G))$ 
such that $\sigma(x) = x$ and for any $f \in {\cal F}_x$, $f$ can 
be transformed into $\sigma \circ f$ by a sequence of finite 
moves. Then (1) $G_x$ is isomorphic to a subgroup of $S_{n-1}$, 
(2) $G_x$ is independent on $x$ up to isomorphism, and 
(3) $c({\rm puz}(G,1)) = [S_{n-1}\,:\,G_x] = (n-1)! / |G_x|$. 

For positive integers $a_1$, $a_2$, $a_3$, we define 
$\theta(a_1, a_2, a_3)$-graph such that 
(1) there exists a pair of vertices $u$ and $v$ of 
degree $3$, and 
(2) $u$ and $v$ are linked by three disjoint paths 
containing $a_1$, $a_2$ and $a_3$ inner vertices, respectively.  
\medskip\\
{\bf Theorem A(Wilson\cite{W1974}).} %
Let $n \ge 2$. 
Let $G$ be a graph with $n$ vertices. 
Suppose that $G$ is $2$-connected and $G$ is not a cycle. 
Let $c = c({\rm puz}(G,1))$. 
\\
{\rm (1)} If $G$ is a bipartite graph, then $G_x \cong A_{n-1}$ 
and $c =2$. 
\\
{\rm (2)} If $G$ is not a bipartite graph except $\theta(1,2,2)$, 
then $G_x \cong S_{n-1}$ and $c = 1$. 
\\
{\rm (3)} If $G$ is $\theta(1,2,2)$, then $G_x \cong PGL_2(5)$ 
and $c=6$, where $PGL_2(5)$ is the projective general linear 
group on $2$-dimensional vector space over a finite field of 
order $5$.
\medskip

Theorem A is generalized for the case $q \ge 2$ \cite{KMS1984}. 
Let $G$ be a connected graph. 
Let $k$ be a positive integer. 
A path $I=v_1 v_2 \ldots v_k$ of $G$ is called an \textit{isthmus} 
if $V(G)\setminus V(I)$ is partitioned into nonempty partite sets 
$X$ and $Y$ such that every path from $X$ to $Y$ passes 
through $I$. In this definition, we say that the isthmus $I$ 
\textit{separates} $X$ and $Y$. An isthmus with $k$ vertices 
is called a $k$-isthmus. Note that a $1$-isthmus is a cut-vertex 
of $G$. 
\medskip\\
{\bf Theorem B(Kornhauser, Miller, Spirakis\cite{KMS1984}).} %
Let $2 \le k \le n-1$. 
Let $G$ be a connected graph with $n$ vertices. 
Suppose that $G$ is not a cycle. 
Then the following conditions are equivalent. 
\\
{\rm (1)} $G$ has no $k$-isthmus. 
\\
{\rm (2)} $(G,k)$ is transitive. 
\\
{\rm (3)} $(G,k)$ is feasible. 
\medskip

In application, it is important to consider the number of 
moves that are necessary, or algorithms to transfer the 
pebbles.  Motion planning on graphs are studied 
intensively\cite{AMPP1999,CDP2006,PRST1994,RW1990}.  

In the next two sections, we will focus on analyzing 
the following three problems for pebble motion on 
graphs, each of which play a key role in the latter 
sections. 
\\
{\rm (1)} \textit{transitivity problem}: What is the reachable set 
of vertices for a given pebble? 
\\
{\rm (2)} \textit{contact problem}: Can a given pair of pebbles 
contact each other? 
\\
{\rm (3)} \textit{equivalence problem}: Can a given pair of 
configurations be equivalent to each other? 

In Section\tume\ref{Isthmus_Structure}, we investigate isthmus 
structure of graphs and introduce the ($k$-)isthmus tree of an 
underlying graph. The isthmus tree describes how the underlying 
graph contains its isthmuses. 
In Section\tume\ref{Contact_Condition}, we will see that the equivalence 
problem\tume{$(3)$} for a given pair of configurations can be reduced 
to the contact problem\tume{$(2)$} for each of the configurations. 
In order to solve these problems, it turns to be useful 
to use the notion of $k$-isthmus tree, where $k$ is the number of 
unoccupied vertices of the board graph. 
Note that an algorithm which solves the equivalence problem $(2)$  
in the above and generates an efficient sequence of moves from one 
configuration to the other, if any, in $\Od(|V(G)|^3)$-time, was 
already anounced by Kornhauser et al.\tume\cite{KMS1984}. 
However, they have never presented its details. 


\section{Isthmus Structure of Graphs}\label{Isthmus_Structure}
Let $G$ be a connected simple graph.
Let us denote a subgraph of $G$ induced by $S \subset V(G)$ by $G[S]$.
A set of vertices $B$ of $G$ is called a $k$-\textit{block} if
(1) $G[B]$ is connected,
(2) $G[B]$ has no $k$-isthmus of $G[B]$, and
(3) $B$ is maximal with respect to (1) and (2).
Namely, for any proper supset $S$ of $B$, $G[S]$ is disconnected or $G[S]$ has a $k$-isthmus of $G[S]$.

Note that a $1$-block is simply a block of a given graph.

We have $|B| \ge k+1$ for any $k$-block $B$ of $G$ with $B \ne V(G)$, because $G[C]$ has no $k$-isthmus for $C \subset V(G)$ with $|C| \le k+1$.
For example, the graph $G$ in Fig.1 has four $3$-isthmuses, and five $3$-blocks.   

\begin{center}
\includegraphics[scale=0.6]{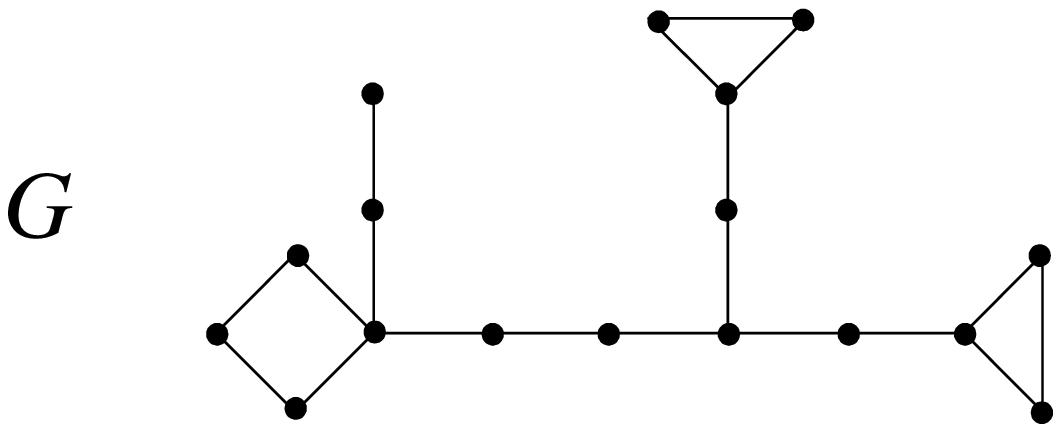}
\vspace*{1cm}\\
\includegraphics[scale=0.6]{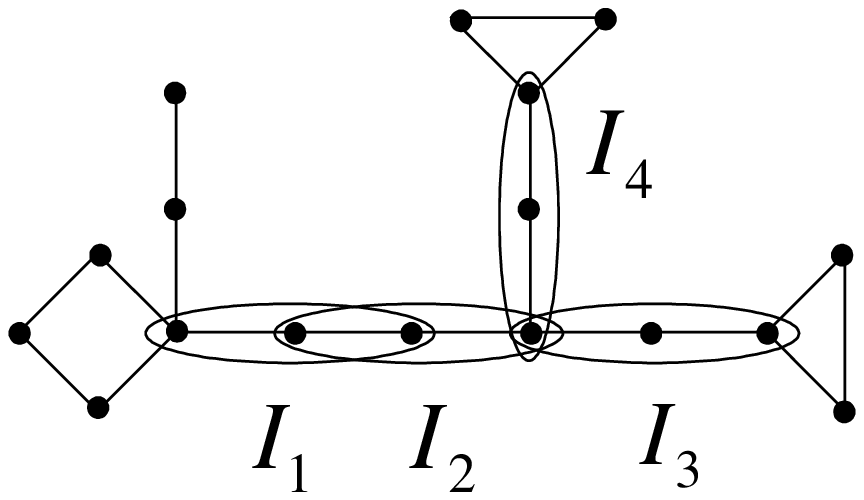}
\vspace*{1cm}\\
\includegraphics[scale=0.6]{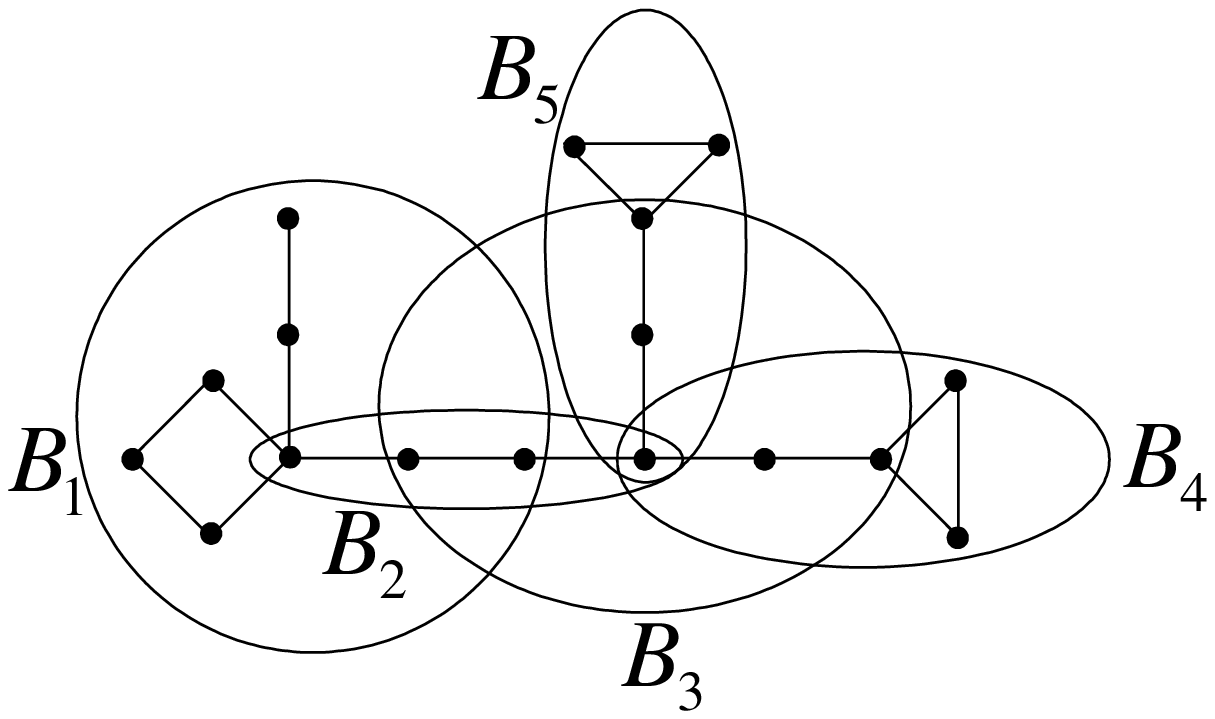}
\vspace*{1cm}\\
\textbf{Fig.1. }{$3$-isthmuses and $3$-blocks of a graph.}
\end{center}
    
\begin{la}\label{k-block}
Let $S \subset V(G)$.
If $|S|=k+1$, then there exists at most one $k$-block $B$ of $G$ such that $S \subset B$.
Moreover, if $|S|=k+1$ and $G[S]$ is connected, then there exists a unique $k$-block $B$ of $G$ such that $S \subset B$.
\end{la}
{\bf Proof.} %
Suppose toward a contradiction that $B_1$ and $B_2$ are two distinct $k$-blocks of $G$ such that $S \subset B_i$ for $i=1,2$.
Let $B=B_1 \cup B_2$.
Since $G[B]$ is connected, by the maximality of $B_2$,
 there exists a $k$-isthmus $I$ of $G[B]$.
Hence, we have a partition $B=X \cup Y \cup V(I)$ such that every path from $X$ to $Y$ is passing through $I$.

Since $|I| < |S|$, there exists a vertex $u \in S \setminus I$.
We may assume $u \in X$.
Take a vertex $v \in Y$.
We may assume $v \in B_1$ without loss of generality.
Then both $u$ and $v$ are contained in $B_1$.
Hence, by the connectivity of $G[B_1]$, we have $I \subset B_1$.
Then $I$ is a $k$-isthmus of $G[B_1]$, which contradicts that $B_1$ is a $k$-block.
Hence, there exists at most one $k$-block $B$ such that $S \subset B$.

Suppose that $G[S]$ is connected.
Since $G[S]$ has no $k$-isthmus, there exists a subset $B \supset S$ such that (1) $G[B]$ is connected, (2) $G[B]$ has no $k$-isthmus of $G[B]$ and (3) $B$ is maximal with respect to (1) and (2).
Then $B$ is a $k$-block satisfying $S \subset B$.
\owari
\medskip

\begin{la}\label{B_covers_S}
Let $S \subset V(G)$ with $|S| \ge k+1$.
If $G[S]$ is connected and $G[S]$ has no $k$-isthmus of $G[S]$,
then there exists a unique $k$-block $B$ such that $B \supset S$.
\end{la}
{\bf Proof.} %
Take a family of $(k+1)$-subsets $\{ S_j \}_{j \in J}$ of $V(G)$ such that $S = \cup_{j \in J} S_j$ and $G[S_j]$ is connected for $j \in J$.
Then, by Lemma 1, there exists a family of $k$-blocks $\{ B_j \}_{j \in J}$ such that $B_j \supset S_j$ for $j \in J$.
Since $S$ has no $k$-isthmus of $G[S]$, we have $B_j \supset S$ for $j \in J$.
Hence, $B_j$ coincides with each other for $j \in J$.
Let $B=B_j$ for $j \in J$.
Then $B$ is a $k$-block containing $S$, as required.
\owari
\medskip

For a vertex $v \in V(G)$, let $N(v)$ denote the set of neighbours of $v$ in $G$.
For a vertex set $S \subset V(G)$, let $N(S) = \cup_{v \in S} N(v) \setminus S$.
For a pair of subsets $T_1$ and $T_2$ of $V(G)$,
a path $P$ of $G$ is called a $(T_1, T_2)$-path if $P$ begins from $T_1$ and ends at $T_2$ and all the inner vertices of $P$ are not contained in $T_1 \cup T_2$.

\begin{la}
Let $B_1$ and $B_2$ be distinct $k$-blocks of a graph $G$.
If there exists a $k$-subset $S \subset B_1 \cap B_2$ such that $G[S]$ is connected,
 then $G[S]$ is a $k$-isthmus of $G[B_1 \cup B_2]$ and $G[S]$ is a $k$-isthmus of $G$.
\end{la}
{\bf Proof.} %
By the maximality of $B_1$, there exists a $k$-isthmus $I = v_1 v_2 \ldots v_k$ of $G[B_1 \cup B_2]$.
We claim that $V(I) = S$.
Suppose to a contradiction that $V(I) \ne S$.
Then there exists a vertex $u \in S \setminus V(I)$,
 because $|S| = |V(I)|$.
Since $I$ is a $k$-isthmus of $G[B_1 \cup B_2]$, 
 there exists a vertex $v \in (B_1 \cup B_2) \setminus I$ such that $u$ and $v$ are separated by $I$.
By symmetry, we may assume $v \in B_1$.
Since $\{ u,v \} \subset B_1$, by the connectivity of $G[B_1]$, we have $I \subset B_1$.
This implies that $G[B_1]$ has a $k$-isthmus of itself, a contradiction.
Therefore, we have $V(I)=S$.

Next, we show that the path $I$ is a $k$-isthmus of $G$.
Since $I$ is a $k$-isthmus of $G[B_1 \cup B_2]$, there is a partition $B_1 \cup B_2 = X \cup Y \cup V(I)$ such that $I$ separates $X$ and $Y$ in $G[B_1 \cup B_2]$.
Let $Z = V(G) \setminus (B_1 \cup B_2)$.
Let $u$, $v$ be two endvertices of $I$ such that $N(u) \cap X \ne \emptyset$ and $N(v) \cap Y \ne \emptyset$.
Let $X_1 = X \cup \{ u \}$, $Y_1 = Y \cup \{ v \}$ and $S_1 = V(I) \setminus \{ u,v \}$.

If there exists a vertex $z \in Z$ such that $z$ is adjacent to $S_1$, then $B_1 \cup \{ z \}$ has no isthmus of itself.
This contradicts to the maximality of $B_1$.
Hence, all $(Z,B_1 \cup B_2)$-paths end at $X_1 \cup Y_1$.
if there exists a vertex $z \in Z$ such that
 there exists a $(z,X_1)$-path and there exists a $(z,Y_1)$-path, then we have a $(X_1, Y_1)$-path $P$ such that $V(P) \cap S_1 = \emptyset$.
Then, $G[B_1 \cup B_2 \cup V(P)]$ has no $k$-isthmus of itself.
This contradicts to the maximality of $B_1$.
Hence, we have a partition $Z = X_2 \cup Y_2$ such that
 any $(X_2,B_1 \cup B_2)$-path ends at $X_1$ and 
 any $(Y_2,B_1 \cup B_2)$-path ends at $Y_1$.
Let $X' = X \cup X_2$ and $Y' = Y \cup Y_2$.
Then we have a partition $V(G)=X' \cup Y' \cup V(I)$ such that $I$ separates $X'$ and $Y'$.
Therefore, $I$ is a $k$-isthmus of $G$.
\owari
\medskip

For a connected graph $G$, let ${\cal I}_k$ and ${\cal B}_k$ denote the set of all $k$-isthmuses of $G$ and the set of all $k$-blocks of $G$.
The $k$-\textit{isthmus graph} $T_k$ of $G$ is defined such that 
$V(T_k)$ is ${\cal B}_k \cup {\cal I}_k$ and 
$E(T_k)$ is $\{ (B,I) \in {\cal B}_k \times {\cal I}_k : V(I) \subset B \}$.

\begin{center}
\includegraphics[scale=0.6]{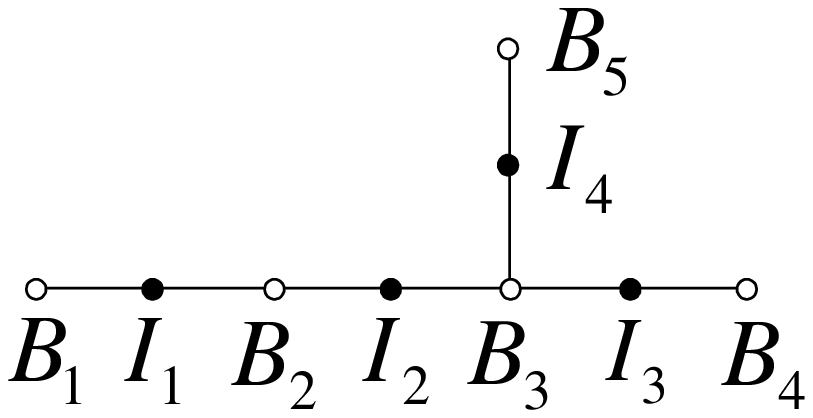}
\\
\textbf{Fig.2. }{The isthmus graph of a graph $G$ in Fig.1.}
\end{center}

By definition, the isthmus graph is a bipartite graph.
As a matter of fact, it is a tree.

\begin{prop}
Let $G$ be a connected graph.
Then the $k$-isthmus graph $T_k$ of $G$ is a tree.
\end{prop}
{\bf Proof.} %
Firstly, we show that $T_k$ is connected.
For a $k$-isthmus $I$ of $G$, we can take a vertex $u \in V(G)\setminus I$ such that $G[I \cup \{ u \}]$ is connected.
By Lemma 1, we have a $k$-block $B$ containing $I$.
Hence, $I$ is not isolated in $T_k$.
 
Let $B_1$ and $B_2$ be distinct $k$-blocks of $G$.
We show that there exists a path from $B_1$ to $B_2$ in $T_k$. 
Take a $(k+1)$-subsets $A_i$ of $B_i$ such that  
 $G[A_i]$ is connected for $i=1,2$.
Since $G$ is connected, we have a finite sequence of $(k+1)$-sets $A_1 = A'_1, A'_2, \ldots, A'_s = A_2$ of $V(G)$ such that
 $G[A'_i]$ is connected for $1 \le i \le s$ and
 $|A'_i \cap A'_{i+1}|=k$, $G[A'_i \cap A'_{i+1}]$ is connected for $1 \le i \le s-1$.

By Lemma 1, we have a sequence of $k$-block $B'_i$ such that $A'_i \subset B'_i$ for $1 \le i \le s$.
Note that $B_1 = B'_1$ and $B_2 = B'_s$.
Furthermore, since $|A'_i \cap A'_{i+1}|=k$ and $G[A'_i \cap A'_{i+1}]$ is connected, by Lemma 3, we have $B'_i = B'_{i+1}$ or $B'_i$ and $B'_{i+1}$ contains a common $k$-isthmus for $1 \le i \le s-1$.
Hence, there is a path from $B_1$ to $B_2$ in $T_k$.
Therefore, $T_k$ is connected. 

Secondly, we will show that $T_k$ has no cycle.
Suppose to a contradiction that $T_k$ has a cycle $C=I_1 B_1 I_2 \ldots I_s B_s I_1$, where $I_i$ is a $k$-isthmus and $B_i$ is a $k$-block for $1 \le i \le s$ such that $B_i \cap B_{i+1} = V(I_{i+1})$ for $1 \le i \le s-1$ and $B_s \cap B_1 = V(I_1)$.
Let us take a partition $V(G)=X \cup Y \cup V(I_1)$ such that $I_1$ separates $X$ and $Y$.
Since $I_1$ is a $k$-isthmus contained in $B_s \cap B_1$, by Lemma 3, we may assume $B_1 \subset X \cup V(I_1)$ and $B_s \subset Y \cup V(I_1)$.
We claim that every $k$-block $B$ of $G$ satisfies $B \subset X \cup V(I_1)$ or $B \subset Y \cup V(I_1)$.
Indeed, otherwise we have $B \cap X \ne \emptyset$ and $B \cap Y \ne \emptyset$.
Since $G[B]$ is connected, we have $B \supset V(I_1)$.
It follows that $I_1$ is a $k$-isthmus of $B$, a contradiction.

Then there exists an index $\alpha$ with $2 \le \alpha \le s$ such that $B_{\alpha-1} \subset X \cup V(I_1)$ and $B_\alpha \subset Y \cup V(I_1)$.
Since $I_\alpha$ is a $k$-isthmus contained in $B_{\alpha-1} \cap B_\alpha \subset V(I_1)$, it coincides with $I_1$, a contradiction.
\owari
\medskip

We call $T_k$ the $k$-\textit{isthmus tree} of a graph $G$.

\section{Contact Condition of Pebbles}\label{Contact_Condition}
Let $G$ be a connected graph with $n$ vertices.
Let $P$ be the set of pebbles of order $m$. 
For the configurations of $P$ on $G$, let us define 
the {\em vacancy size}, denoted by $k(P,G):=n-m$,  as 
the number of unoccupied vertices of $G$. In this section, 
let us use $k$ in stead of $k(P,G)$, as an abbreviation. 
Because we treat only configurations which admit a move 
of a pebble, let us assume $k \geqq 1$ throughout this section. 

For a configuration $f \in \cal{F}$ and a pebble $p \in P$,
let us define $R(p,f) = \{ v \in V(G) : g(p)=v \textrm{~for some
  configuration $g$ such that~} g \sim f \}$.
We call $R(p,f)$ the \textit{reachable range} of $p$ starting from $f$.

For a configuration $f \in \cal{F}$ and a pebble $p \in P$, 
let $v = f(p)$. Since $G$ is connected, we can gather $k$  
unoccupied vertices around $v$ without moving $p$.  
More precisely, we have a configuration $f_1$ equivalent to  
$f$ such that (1) $f_1(p)=v$, (2) $G[A]$ is connected where 
$A = V(G) \setminus f_{1}(P\setminus\{p\})$. 
Note that $A$ is not uniquely determined by a given 
pair $p \in P$ and $f \in \cal{F}$. 
Since $|A|=k+1$, by Lemma 1 in the previous section, 
we have a $k$-block $B$ of $G$ containing $A$. 

\begin{thm}\label{Range=Block}
Let $p \in P$ and $f \in \cal{F}$.
Let $B$ be a $k$-block as defined as above.
Then $R(p,f)=B$.
\end{thm}
{\bf Proof.} %
Firstly, we will show that $B \subset R(p,f)$.
Since $G[B]$ is connected and $G[B]$ has no 
$k$-isthmus of itself, by Theorem A and Theorem B, 
 $p$ can be moved to all the vertices of $B$. 
Hence, we have $B \subset R(p,f)$. 

Secondly, we will show that $R(p,f) \subset B$. Suppose 
to a contradiction that $R(p,f) \setminus B \ne \emptyset$. 
Since $G[R(p,f)]$ is connected, we may assume there exists 
a vertex $u \in R(p,f) \setminus B$ such that 
$G[B \cup \{ u \}]$ is connected. Then by the maximality of 
$B$, there exists a $k$-isthmus $I=v_1 v_2 \ldots v_k$ of 
$G[B \cup \{ u \} ]$. Then we have $V(I) \subset B$.
Let $A=V(I) \cup \{ u \}$. 
Since $G[A]$ is connected and $|A|=k+1$, 
by Lemma\tume\ref{k-block}, 
we have a $k$-block $B'$ such that $V(I) \cup \{ u \} \subset B'$. 
Then, by Lemma 3, $I$ is a $k$-isthmus of $G$. 
We may assume $u$ is a neighbour of $v_k$. 
Let $C$ be the set of vertices of $G$ separated by $I$ from $u$. 
Hence, for any configuration $g$ equivalent to $f$ such that 
$g(p)=v_i$ for some $1 \le i \le k$, the number of unoccupied 
vertices of $C \cup \{ v_{1}, v_{2}, \ldots, v_{i-1} \}$ is at least
$i$. In particular, if $g \sim f$ and $g(p)=v_k$, then all the 
unoccupied vertices are contained in $C \cup I$.  
Therefore, $p$ cannot be moved to $u$, a contradiction.
\owari
\medskip

Next, we consider contact condition. 
Let $p$, $q$ be a pair of distinct pebbles. 
We say that $p$ \textit{contacts} $q$ beginning from an initial 
configuration $f$, if there exists a configuration $g$ 
equivalent to $f$ such that $g(p)$ and $g(q)$ are adjacent.

If $G$ is a cycle, it is easy to see that $p$ contacts $q$ 
if and only if $q$ is next to $p$ along the cycle in the 
initial configuration. If $k \ge 2$ and a $k$-block $B$ 
is a proper subset of $V(G)$, then $G[B]$ is not a cycle. 

\begin{thm}\label{contact_decision}
Let $p$, $q \in P$, and $f \in \cal{F}$.
Then $p$ can contact $q$, if and only if one of 
the following conditions hold;
 (1) $R(p,f)=R(q,f)$ and $G[R(p,f)]$ is not a cycle, or 
 (2) $R(p,f)=R(q,f)$ and $G[R(p,f)]$ is a cycle and $q$ is 
next to $p$ along the cycle, or
 (3) $R(p,f) \ne R(q,f)$ and $G[R(p,f) \cap R(q,f)]$ is 
a $k$-isthmus of $G$.
\end{thm}
{\bf Proof.} %
Firstly, we show that any of the conditions (1), (2) or (3) is
sufficient for the contact. 
Suppose that $R(p,f)=R(q,f)$.
Let $B=R(p,f)$.
By Theorem\tume\ref{Range=Block}, 
$B$ is a $k$-block of $G$.
Then we have a configuration $f_1$ equivalent to $f$ 
such that $V(G) \setminus f_1(P\setminus\{p\}) \subset B$.
Since $R(p,f)=R(q,f)$, we have 
$V(G) \setminus f_1(P\setminus\{p,q\}) \subset B$.
Let us consider the puzzle on the restricted graph $G[B]$.
Let $m'$ be the number of pebbles on $B$ in $f_1$. 
It is easy to see that if $G[B]$ is a cycle and $q$ is next to 
$p$ along the cycle, $p$ can contact $q$. 
Suppose that $G[B]$ is not a cycle. 
By Theorem A and Theorem B, the puzzle 
$\textrm{puz}(G[B],k)$ is $3$-transitive for $m' \ge 3$. 
For $m'=2$, since $G[B]$ is not a cycle, 
$\textrm{puz}(G[B],k)$ is feasible.
Hence, in all the cases, $p$ can contact $q$.

Next, we assume that $R(p,f) \ne R(q,f)$ and 
$I=G[R(p,f) \cap R(q,f)]$ is a $k$-isthmus $I$ of $G$. 
Let us take two vertices $u \in (R(p,f) \setminus I) \cap N(I)$ 
and $v \in (R(q,f) \setminus I) \cap N(I)$. 
Since $\textrm{puz}(G[R(p,f)],k)$ is feasible, there 
exists a configuration $f_1$ equivalent to $f$ such that 
$f_1(p)=u$, and $V(G) \setminus f_1(P)=I$.
Then since the puzzle $\textrm{puz}(G[R(q,f_1)],k)$ is 
feasible, there exists a configuration $f_2$ equivalent 
to $f_1$ such that $f_2(p)=u$, $f_2(q)=v$, and 
$V(G)\setminus f_2(P)=I$.
Then by using a path $I$, $p$ can contact $q$.

Secondly, we show that one of the conditions (1), (2) or (3) 
is necessary for the contact. Let $f_1$ be a configuration 
equivalent to $f$ such that $u=f_1(p)$ and $v=f_1(q)$, 
where $u$ and $v$ are adjacent. Let us gather $k$ 
unoccupied vertices around $u, v$ without moving 
$p$ and $q$. More precisely, we have a configuration 
$f_2$ equivalent to $f_1$ such that (1) $f_2(p)=u$, 
(2) $f_2(q)=v$, (3) $G[A]$ is connected where 
$A = V(G) \setminus f_2(P \setminus \{p,q \})$.
If $G[A]$ has no $k$-isthmus of $G$, 
by Lemma\tume\ref{B_covers_S}, 
we have a unique $k$-block $B$ containing $A$. 
In this case, by Theorem 5, we have $B=R(p,f)=R(q,f)$. 
If $G[R(p,f)]$ is a cycle and $q$ is not next to $p$ 
along the cycle, $p$ cannot contact $q$. 

Suppose that $G[A]$ has a $k$-isthmus $I$ of $G$.
In this case, $G[A]$ is a path of $G$.
Let $x_1$, $x_2$ be the two endvertices of $G[A]$.
Note that $A \setminus \{ x_1, x_2 \} = V(I)$.
Then there exists a configuration $f_3$ equivalent to 
$f_2$ such that $f_3(\{p,q\} )=\{ x_1, x_2 \}$ and 
$f_3^{-1}(I)=\emptyset$.
Then we have $G[R(p,f) \cap R(q,f)] = I$, as required.
\owari
\medskip

Lastly, we consider a necessary and sufficient 
condition such that the two configurations $f$ 
and $g$ are equivalent. 

First, we deal with the case $k \ge 2$.

\begin{thm}\label{equivalence_decision_k_ge_2}
Let $k \ge 2$.
Let $f$ and $g$ be two configurations.
Then $f$ and $g$ are equivalent if and only if
all the following conditions hold; 

\begin{u-enumerate}
\item{$R(p,f)=R(p,g)$ for any pebble $p \in P$. }
\item{If $G$ is a cycle graph, then the cyclic order of $P$ 
on $G$ is the same as in $f$ and $g$. }
\end{u-enumerate}
\end{thm}
{\bf Proof.} %
Suppose that $f$ and $g$ are equivalent.
Then any pebble $p$ on $f(p)$ with a configuration 
$f$ can be moved at $g(p)$ with $g$, and vice versa.
Hence we have $R(p,f) = R(p,g)$.
It is not difficult to check the condition (2).

Conversely, suppose that the conditions (1) and (2) hold.
We proceed by induction on the number $s$ of $k$-blocks of $G$.
Let $s=1$.
Since $G$ has no $k$-isthmus, by Theorem B, if $G$ is not a cycle, $(G,k)$ is feasible.
Hence, $f$ and $g$ are equivalent.
If $G$ is a cycle, by the condition (2), $f$ and $g$ are equivalent.

Let $s \ge 2$.
In this case, note that any $k$-block of $G$ is not a cycle, because it contains a set of vertices of some $k$-isthmus.
Let us take a $k$-block $B$ such that $B$ corresponds to a leaf of the $k$-isthmus tree of $G$.
Let $I$ be a unique $k$-isthmus such that $B \supset V(I)$.
Then we have two configurations $f_1$ and $g_1$ such that $f_1 \sim f$, $g_1 \sim g$ and $f_1^{-1}(V(I))=g_1^{-1}(V(I))=\emptyset$.
Put $P_1 = f_1^{-1}(B)$.
Then, for $p \in P$, we have $p \in P_1$ if and only if $R(p,f)=B$.
By the condition (1), we have $P_1 = g_1^{-1}(B)$.
Since $(G[B].k)$ is feasible, $f_1|_{P_1}$ and $g_1|_{P_1}$ are equivalent.
Now, we have two configurations $f_2$ and $g_2$ such that $f_2 \sim f$, $g_2 \sim g$, and $f_2(p)=g_2(p) \in B \setminus V(I)$ for all $p \in P_1$.
Let us define $G'=G \setminus (B \setminus V(I))$ and $P'=P \setminus P_1$.
Then $G'$ has $s-1$ $k$-isthmuses.
By inductive hypothesis, $f_2|_{P'}$ and $g_2|_{P'}$ are equivalent on $G'$.
Therefore, $f$ and $g$ are equivalent on $G$.
\owari
\medskip

Next, we deal with the case $k = 1$.
In this case, according to Theorem A, the conditions for equivalence becomes a little complicated.

\begin{thm}\label{equivalence_decision_k_eq_1}
Let $k = 1$.
Let $f$ and $g$ be two configurations.
Then $f$ and $g$ are equivalent if and only if
all the following conditions hold; 

\begin{u-enumerate}
\item{$R(p,f)=R(p,g)$ for any pebble $p \in P$. }
\item{For $x \in V(G)$, 
let $f_x$ and $g_x$ be arbitrary configurations such that 
$f_x \sim f$, $g_x \sim g$ and $f_x^{-1}(x)= g_x^{-1}(x)= \emptyset$. 
\begin{u-enumerate}
    \item{Suppose that $G[R(p,f)]$ is a cycle graph for a pebble $p$.
Let $x \in R(p,f)$.
Then $f_x^{-1}(y)=g_x^{-1}(y)$ for all $y \in R(p,f) \setminus \{ x \}$. }
    \item{Suppose that $G[R(p,f)]$ is a bipartite graph for a pebble $p$.
Let $x \in R(p,f)$.
Then, $g_x\circ f_x^{-1}$ restricted 
on $R(p,f) \setminus \{ x \}$ is an even permutation. }
    \item{Suppose that $G[R(p,f)]$ is the $\theta(1,2,2)$ graph for a pebble $p$.
Let $x \in R(p,f)$.
Then, $g_x\circ f_x^{-1}$ restricted 
on $R(p,f) \setminus \{ x \}$ is contained in $PGL_2(5)$, which is the projective general linear group on $2$-dimensional 
vector space over a finite field of order $5$. }
\end{u-enumerate}
}
\end{u-enumerate}
\end{thm}

Theorem 8 is proved in a similar manner as in the proof of Theorem 7, and we omit the proof.


\section{Equivalence of Arrangements}\label{equivalency}

In this section, we focus on the case of AAP. The proof for the case 
of SGA is almost identical to the case of AAP and will be omitted.  

\begin{la}\label{contraction_of_agents}
Let $(G_A, G_M)$ be an ordered pair of graphs, let $\{ps,t\}$ be an 
edge of $G_M$, and let $f$ be an arrangement of $(G_A, G_M)$ 
such that \mbox{$|f^{-1}(\{s,t\})| \geqq 2$} and the graph
$G_A[f^{-1}(\{s,t\})]$ is connected. Then the arrangement $f$ is 
t-equivalent (in the sense of AAP) to the following arrangement $g$:
\begin{equation*}
g(x)= \begin{cases}
            t, & \text{if $f(x) \in \{s,t\}$;} \\
            f(x), & \text{if $f(x) \notin \{s,t\}$.}
            \end{cases}
\end{equation*}
\end{la}
Note that, as for the case of SGA, by the definition, 
the conclusion of this lemma holds even if we omit 
the assumption ``$|f^{-1}(\{s,t\})| \geqq 2$'' for the 
arrangement $f$ from its statement. 
\begin{Proof}
We divide our proof into two cases.
\begin{case}{\em $f^{-1}(t) \neq \emptyset$.

\noindent
Let $U:=f^{-1}(s)$. Because $f^{-1}(\{s,t\}) \setminus U = f^{-1}(t)
\neq \emptyset$ and $G_A[U]$ is connected, by the definition 
of transfer of AAP, the arrangement $f$ is t-equivalent to the 
arrangement $g$. 
}\end{case}

\begin{case}{\em $f^{-1}(t)=\emptyset$.

\noindent
In this case, because $f^{-1}(s)\geqq 2$ and $G_A[f^{-1}(s)]$ 
is connected, there exists at least one vertex $u$ in $f^{-1}(s)$ 
such that $G_A[f^{-1}(s)]-u$ remains connected. Thus, by 
the definition of transfer of AAP, the arrangement $f$ is 
t-equivalent to the following arrangement $f^{\prime}:$ 
\begin{equation*}
f^{\prime}(x)= \begin{cases}
            t, & \text{if $x=u$;} \\
            f(x), & \text{if $x \neq u$.}
            \end{cases}
\end{equation*}
This $f^{\prime}$ satisfies the condition of {\bf Case 1}, 
and hence $g \cong f^{\prime} \cong f$ holds. 
}\end{case}
\end{Proof}

Let $f$ be an arrangement of an ordered pair $(G_A, G_M)$ of graphs. 
Lemma\tume\ref{contraction_of_agents} indicates that, if a country 
$s$ of the route map $G_M$ contains at least two agents under 
the arrangement $f$, we can treat the agent sub-network 
$G_A[f^{-1}(s)]$ as if a single `pebble' on the vertex $s$ of 
the `board graph' $G_M$. This observation leads us to define 
the following set of new concepts for the AAP model. For every 
arrangement $f$ of the ordered pair $(G_A, G_M)$, let us define 
the {\em isolated-agent set} $\IA_{f}$ as the set of agents
$\big\{a \in V(G_A) : f^{-1}(f(a))=\{a\}\big\}$. In other words, 
the set $\IA_{f}$ is the set of all agents who are staying 
alone in their countries under the arrangement $f$.  In the same 
way, let us define the {\em isolated-country set} $\IC_{f}$ as 
the set of countries $\big\{c \in V(G_M) : |f^{-1}(c)|=1\big\}$. 
Note that $f(\IA_{f})=\IC_{f}$ holds. 
Now, let us consider the following configuration $\phi_f$ of the 
pebble motion problem corresponding to the arrangement $f$: 
The pebble set $P_{\phi_f}$ for $\phi_f$ is the set 
$\{f^{-1}(c) : c \in V(G_M)\setminus \IC_{f}, f^{-1}(c) \neq 
\emptyset \}$. The board graph $G_{\phi_f}$ for $\phi_f$ is 
the subgraph of $G_M$ induced by the vertex set 
$V(G_M)\setminus \IC_{f}$. Then our configuration is 
defined by $\phi_f : P_{\phi_f} \ni f^{-1}(c) \mapsto c \in V(G_{\phi_f})$. 
Let us call this $\phi_f$ the {\em configuration associated 
with $f$}. For two arrangements $f$ and $g$ of $(G_A, G_M)$,  
we say that the configuration $\phi_f$ associated with $f$ is 
{\em equivalent}  to the configuration $\phi_g$ associated 
with $g$ and use the notation $\phi_f \sim \phi_g$, if and 
only if both their pebble sets and board graphs are coincident 
with each other ($P_{\phi_f}=P_{\phi_g}$ and $G_{\phi_f}=G_{\phi_g}$) 
and $\phi_f$ is equivalent to $\phi_g$ in the sense of the pebble 
motion problem. It is clear that $\phi_f \sim \phi_g$ implies 
$f \cong g$. Although the converse is not true in general, 
the following extremal condition warrants its affirmation: 

Let $f$ be an arrangement of an ordered pair $(G_A, G_M)$ 
of graphs. Let $\phi_f$ denote the configuration associated with $f$, 
$P_{\phi_f}$ its pebble set, $G_{\phi_f}$ its board graph. Now 
let us define the set of {\em $f$-irreducible arrangements} $\Irr(f):=
\big\{ h : h \cong f, |P_{\phi_h}|=\min\{ |P_{\phi_g}| : g \cong f\}\big\}$ 
and the set of {\em $f$-irreducible configurations} $\Phi(\Irr(f)):=
\{ \phi_h : h \in \Irr(f)\}$. 
\begin{thm}\label{minimum_contraction}
For an arbitrary arrangement $f$ of an arbitrary ordered pair 
$(G_A, G_M)$ of graphs, all configurations in the set $\Phi(\Irr(f))$ 
are equivalent. 
\end{thm}

In order to prove the above theorem, we shall prepare some 
notations and prove a technical theorem. 

Let $(G_A, G_M)$ be an ordered pair of graphs, $\{s,t\}$ be an 
edge of $G_M$, and let $f$ be an arrangement of $(G_A, G_M)$ 
such that \mbox{$f^{-1}(s) \neq \emptyset, |f^{-1}(\{s,t\})| \geqq 2$} 
and the graph $G_A[f^{-1}(\{s,t\})]$ is connected. 
Lemma\tume\ref{contraction_of_agents} guarantees 
that the following arrangement $g$ can be achieved from 
the arrangement $f$ by at most two steps of transfers 
(in the sense of AAP):
\begin{equation*}
g(x)= \begin{cases}
            t, & \text{if $f(x) \in \{s,t\}$;} \\
            f(x), & \text{if $f(x) \notin \{s,t\}$.}
            \end{cases}
\end{equation*}
Let us call the sequence of (at most two) transfers from the
arrangement $f$ to the arrangement $g$ a {\em SUBNET 
MERGER}. Furthermore, let us call a SUBNET MERGER from 
the arrangement $f$ to the arrangement $g$ a {\em SUBNET 
MOVE} if the set $f^{-1}(t)$ is empty. A SUBNET MERGER is called 
{\em proper} if it is not a SUBNET MOVE. 

\begin{thm}\label{only_MERGER}
Let $f$ be an arbitrary arrangement of an arbitrary ordered pair 
$(G_A, G_M)$ of graphs. Then, for every $f$-irreducible arrangement 
$g$, there exists a sequence $f=:f_0 \cong f_1 \cong \cdots 
\cong f_k:=g$ of t-equivalent arrangements on $(G_A, G_M)$
such that, for all $i \in \{0,\ldots, k-1\}$, each sequence of the 
transfers from the arrangement $f_i$ to the arrangement $f_{i+1}$ 
is a SUBNET MERGER. 
\end{thm}

\begin{Proofof}{Theorem\tume\ref{only_MERGER}}
Because $g$ is an $f$-irreducible arrangement, there exists 
a sequence $f=:f_0 \cong f_1 \cong \cdots \cong f_k:=g$ 
of t-equivalent arrangements such that, for all $i \in \{0,\ldots,
k-1\}$, the sequence of the transfers from $f_i$ to $f_{i+1}$ is either 
a SUBNET MERGER or consisting of a single transfer which is not a 
SUBNET MERGER. For the above sequence $f_0 \cong f_1 \cong 
\cdots \cong f_k$, let $m$ be the minimum number in $\{0,\ldots, k\}$ 
on condition that the sequence of the transfers from $f_{m}$ to 
$f_{k}(=g)$ can be written as a sequence of only SUBNET MERGERs. 
We will show that $m=0$ by reductio ad absurdum. 

\noindent
Because of the minimality of $m$, we can assume that each arrangement 
$f_{i+1} ~(i=m, \ldots,k-1)$ can be achieved from the arrangement 
$f_{i}$ by a SUBNET MERGER. And hence we have that; 
\begin{equation*}
\forall i \in \{m,\ldots,k-1\}, \exists \{s_i,t_i\}\in E(G_M), 
f_{i+1}(x)= \begin{cases}
            t_i, & \text{if $f_{i}(x) \in \{s_i, t_i\}$;} \\
            f_{i}(x), & \text{if $f_{i}(x) \notin \{s_i, t_i\}$.}
            \end{cases}
\end{equation*}
Now suppose that $m \geqq 1$. Since $m$ is minimum, we have that 
the sequence of the transfers from the arrangement $f_{m-1}$ to the 
arrangement $f_{m}$ is consisting of a single transfer which is not a 
SUBNET MERGER. That is, there exist an edge $\{s_{m-1}, t_{m-1}\}$ of 
$G_M$ and a proper non-empty subset $U_{m-1}$ of the set 
${f_{m-1}}^{-1}(s_{m-1})$ (i.e. $\emptyset \neq U_{m-1} \subsetneq
{f_{m-1}}^{-1}(s_{m-1})$) such that: 
\begin{equation*}
f_{m}(x)= \begin{cases}
            t_{m-1}, & \text{if $x \in U_{m-1}$;} \\
            f_{m-1}(x), & \text{otherwise.}
            \end{cases}
\end{equation*}
Now let $h_1$ be the following arrangement: 
\begin{equation*}
h_{1}(x):= \begin{cases}
            t_{m-1}, & \text{if $f_{m-1}(x) \in \{s_{m-1}, t_{m-1}\}$;} \\
            f_{m-1}(x), & \text{if $f_{m-1}(x) \notin \{s_{m-1}, t_{m-1}\}$.}
            \end{cases}
\end{equation*}
Clearly, this $h_{1}$ can be achieved from the arrangement $f_{m-1}$ 
by a SUBNET MERGER. Let $\cal{H}$ be the set of all arrangements 
achieved from $h_1$ by a sequence of only SUBNET MERGERs. 
And let $g^{\prime}$ be an arrangement in $\cal{H}$ such that 
$|g^{\prime}(V(G_A))| = \min\{|h(V(G_A))| :  h \in \cal{H} \}$ holds.
Let $f_{m-1}=:h_0 \cong h_1 \cong \cdots \cong h_l:=g^{\prime}$ be 
a sequence of arrangements such that, for all $i \in \{1, \ldots, l\}$, 
$h_i$ can be achieved from $h_{i-1}$ by a single SUBNET MERGER. 
Then we have that; 
\begin{equation*}
\forall i \in \{0,\ldots,l-1\}, \exists \{u_i,v_i\}\in E(G_M), 
h_{i+1}(x)= \begin{cases}
            v_i, & \text{if $h_{i}(x) \in \{u_i, v_i\}$;} \\
            h_{i}(x), & \text{if $h_{i}(x) \notin \{u_i, v_i\}$.}
            \end{cases}
\end{equation*}
Because of the minimality of $m$, we have that $g \neq g^{\prime}$. 
Furthermore, since $g$ is an $f$-irreducible arrangement, 
$|g^{\prime}(V(G_A))| \geqq |g(V(G_A))|$ holds. And hence there exist 
two distinct agents $a, b ~ (\in V(G_A))$ such that both $g(a)=g(b)$ 
and $g^{\prime}(a) \neq g^{\prime}(b)$ hold. Now let us number all the 
vertices of $G_A$ so that $V(G_A):=\{{a_1:=a}, {a_2:=b},$ $\ldots, 
a_n\}$. Let $\id(S):=\min\{i : a_i \in S\}$ be the function from the 
power set of $V(G_A)$ to the set $\{1,\ldots,n\}$ which returns the 
minimum index number of vertices in a given subset $S$ of $V(G_A)$. 
Corresponding to the above sequence of arrangements $f_{m-1}=h_0 
\cong h_1 \cong \cdots \cong h_l=g^{\prime}$, we will define another 
sequence of arrangements $g^{\prime}=:h^{\prime}_{l} \cong
h^{\prime}_{l-1} \cong \cdots \cong h^{\prime}_{0}$ such that, for all 
$i \in \{0,\ldots,l-1\}$,  
\begin{equation*}
h^{\prime}_{i}(x):= \begin{cases}
            u_{i}, & \text{if $h^{\prime}_{i+1}(x)=v_{i}$ and 
                         $h_{i}(a_{\id({h^{\prime}_{i+1}}^{-1}(v_i))})=u_{i}$;} \\
            h^{\prime}_{i+1}(x), & \text{otherwise.}
            \end{cases}
\end{equation*}
Because each arrangement $h_{i+1} ~(i=0, \ldots,l-1)$ is achieved from 
$h_{i}$ by a SUBNET MERGER, the corresponding sequence of the 
transfers from $h^{\prime}_{i+1}$ to $h^{\prime}_{i}$ defined above is a 
SUBNET MOVE. And hence, for all $i \in \{0,\ldots,l\}$, 
$h^{\prime}_{i} \in \cal{H}$ and $|h^{\prime}_{i}(V(G_A))| =
|g^{\prime}(V(G_A))|$. Here we note that, there exists another representation 
of the arrangements $h^{\prime}_{i} ~(i=0, \ldots,l)$, as follows: 
$$\forall (i,x) \in \{0,\ldots,l\} \times V(G_A), 
h^{\prime}_{i}(x)= h_{i}(a_{\id({g^{\prime}}^{-1}(g^{\prime}(x)))}).$$
Let $f^{\prime}_{m}$ be the following arrangement: 
\begin{equation*}
f^{\prime}_{m}(x):= \begin{cases}
            t_{m-1}, & \text{if $h^{\prime}_0(x)=s_{m-1}$ and 
                         $f_{m}(a_{\id({h^{\prime}_0}^{-1}(s_{m-1}))})=t_{m-1}$;} \\
            h^{\prime}_{0}(x), & \text{otherwise.}
            \end{cases}
\end{equation*}
Because $h^{\prime}_{0}$ is achieved from $f_{m-1} ~ (=h_0)$ 
by a sequence of only SUBNET MERGERs, for every country 
$c \in V(G_M)$, if ${f_{m-1}}^{-1}(c) \neq \emptyset$ then 
${f_{m-1}}^{-1}(c) \subseteq {h^{\prime}_{0}}^{-1}(c)$ holds. 
In particular, if ${f_{m-1}}^{-1}(s_{m-1}) \neq \emptyset \neq 
{f_{m-1}}^{-1}(t_{m-1})$, then ${f_{m-1}}^{-1}(s_{m-1}) \subseteq
{h^{\prime}_{0}}^{-1}(s_{m-1})$ and ${f_{m-1}}^{-1}(t_{m-1}) 
\subseteq {h^{\prime}_{0}}^{-1}(t_{m-1})$ hold, and the graph 
$G_A[{h^{\prime}_{0}}^{-1}(\{s_{m-1},t_{m-1}\})]$ is connected. 
Then $f^{\prime}_{m}$ is achieved from $h^{\prime}_{0}$ 
by a SUBNET MERGER and $|f^{\prime}_{m}(V(G_A))| =
|h^{\prime}_{0}(V(G_A))| -1 = |g^{\prime}(V(G_A))| -1$ holds, which 
contradicts the minimality of the size $|g^{\prime}(V(G_A))|$. 
Hence we have that at least one of ${f_{m-1}}^{-1}(s_{m-1})$ or 
${f_{m-1}}^{-1}(t_{m-1})$ is empty, and that $f^{\prime}_{m}$ 
is achieved from $h^{\prime}_{0}$ by a SUBNET MOVE.
Furthermore, $f^{\prime}_{m}$ turns to be achieved from 
$g^{\prime}$ by a sequence of only SUBNET MOVEs.
Now, let $f^{\prime}_{m-1}:=h^{\prime}_0$ and 
we will define one more sequence of arrangements 
$h^{\prime}_{0}=f^{\prime}_{m-1} \cong
f^{\prime}_{m} \cong \cdots \cong f^{\prime}_{k}$ such that, 
for all $i \in \{m,\ldots,k\}$,  
\begin{equation*}
f^{\prime}_{i}(x):= \begin{cases}
            t_{i-1}, & \text{if $f^{\prime}_{i-1}(x)=s_{i-1}$ and 
                         $f_{i}(a_{\id({f^{\prime}_{i-1}}^{-1}(s_{i-1}))})=t_{i-1}$;} \\
            f^{\prime}_{i-1}(x), & \text{otherwise.}
            \end{cases}
\end{equation*}
Again, we have another representation of the arrangements 
$f^{\prime}_{i} ~(i=m, \ldots,k)$, as follows: 
$$\forall (i,x) \in \{m,\ldots,k-1\} \times V(G_A), 
f^{\prime}_{i}(x)= f_{i}(a_{\id({g^{\prime}}^{-1}(g^{\prime}(x)))}).$$
Because $f^{\prime}_{m}$ is achieved from $f_{m-1}$ by a 
sequence of only SUBNET MERGERs, for every country 
$c \in V(G_M)$, if ${f_{m-1}}^{-1}(c) \neq \emptyset$ then 
${f_{m-1}}^{-1}(c) \subseteq {f^{\prime}_{m}}^{-1}(c)$ holds.
Then because each arrangement $f_{i+1} ~(i=m, \ldots,k-1)$ is 
achieved from $f_{i}$ by a SUBNET MERGER, the corresponding 
sequence of the transfers from $f^{\prime}_{i}$ to $f^{\prime}_{i+1}$ 
defined above is also a SUBNET MERGER. Combining this fact 
with the minimality of the size $|g^{\prime}(V(G_A))|$, we have that 
each sequence of the transfers from $f^{\prime}_{i} ~(i=m, \ldots,k-1)$ 
to $f^{\prime}_{i+1}$ defined above is not only a SUBNET MERGER 
but also a SUBNET MOVE. Then $f^{\prime}_{k}$ turns out to be 
achieved from $g^{\prime}$ by a sequence of only SUBNET 
MOVEs, and hence $f^{\prime}_{k}(a) \neq f^{\prime}_{k}(b)$. 
However it contradicts the fact that 
\begin{align*}
f^{\prime}_{k}(a)&= f_{k}(a_{\id({g^{\prime}}^{-1}(g^{\prime}(a)))})
= f_{k}(a_{\id({g^{\prime}}^{-1}(g^{\prime}(a_1)))})\\
&= f_{k}(a_{1})=g(a_1)=g(a)=g(b)=g(a_2)=f_{k}(a_2)\\
&= f_{k}(a_{\id({g^{\prime}}^{-1}(g^{\prime}(a_2)))})
= f_{k}(a_{\id({g^{\prime}}^{-1}(g^{\prime}(b)))})\\
&=f^{\prime}_{k}(b).
\end{align*}
Hence we have that $g=g^{\prime}$ and $m=0$, which 
completes the proof. 
\end{Proofof}

\begin{Proofof}{Theorem\tume\ref{minimum_contraction}}
Let $g,h$ be two arbitrary $f$-irreducible arrangements, 
$\phi_g, \phi_h$ their corresponding $f$-irreducible configurations. 
We will show that $\phi_g \sim \phi_h$.

\noindent
Because $g$ and $h$ are $f$-irreducible arrangements, from 
Theorem\tume\ref{only_MERGER}, we have that there exists two 
sequences of arrangements $f=:g_0 \cong g_1 \cong \cdots 
\cong g_k:=g$ and $f=:h_0 \cong h_1 \cong \cdots \cong h_l:=h$ 
such that each $g_{i+1} ~(i=0, \ldots, k-1)$  (resp. $h_{j+1} ~(j=0, \ldots, l-1)$ )
is achieved from $g_{i}$ (resp. $h_{j}$) by a single SUBNET MERGER. 
In the same way of the proof of Theorem\tume\ref{only_MERGER}, 
let us number all the vertices of $G_A$ as $V(G_A):=\{a_1,\ldots,a_n\}$ 
and define the function $\id(S):=\min\{i : a_i \in S\}$. 
Corresponding to the above sequence of arrangements
$f=g_0 \cong g_1 \cong \cdots \cong g_k=g$,  
we will define another new sequence of arrangements  
$g=:g^{\prime}_k \cong g^{\prime}_{k-1} \cong \cdots \cong g^{\prime}_0$ 
as follows: 
$$\forall (i,x) \in \{0,\ldots,k\} \times V(G_A), ~ 
g^{\prime}_{i}(x):= g_{i}(a_{\id({g^{\prime}}^{-1}(g^{\prime}(x)))}).$$
Because each arrangement $g_{i+1} ~(i=0, \ldots,k-1)$ is achieved from 
$g_{i}$ by a SUBNET MERGER, the corresponding sequence of the 
transfers from $g^{\prime}_{i+1}$ to $g^{\prime}_{i}$ defined above is a 
SUBNET MOVE. 
Then next, corresponding to the above sequence of arrangements
$f=h_0 \cong h_1 \cong \cdots \cong h_l=h$,  
we will define the following new sequence of arrangements  
$g^{\prime}_{0}=:h^{\prime}_0 \cong h^{\prime}_{1} \cong \cdots \cong h^{\prime}_l$: 
$$\forall (j,x) \in \{0,\ldots,l\} \times V(G_A), ~
h^{\prime}_{j}(x):= h_{j}(a_{\id({g^{\prime}}^{-1}(g^{\prime}(x)))}).$$
Again, because the arrangement $g^{\prime}_{0} ~ (=h^{\prime}_{0})$ 
is achieved from the arrangement $g_k ~ (=g)$ by a sequence of 
only SUBNET MOVEs, and because each arrangement $h_{i+1} ~
(i=0, \ldots,l-1)$ is achieved from $h_{i}$ by a SUBNET MERGER, 
the corresponding sequence of the transfers from $h^{\prime}_{i+1}$ 
to $h^{\prime}_{i}$ defined above is a SUBNET MERGER. Furthermore, 
because $g$ is an $f$-irreducible arrangement, the sequence of the 
transfers from $h^{\prime}_{i+1}$ to $h^{\prime}_{i}$ is not only a
SUBNET MERGER, but also a SUBNET MOVE. And hence 
the arrangement $h^{\prime}_{l}$ is achieved from the arrangement 
$g$ by a sequence of only SUBNET MOVEs. This fact tells us that, 
for every country $c \in V(G_M)$, if $g^{-1}(c) \neq \emptyset$ 
then $g^{-1}(c) = {h^{\prime}_{l}}^{-1}(h(a_{\id(g^{-1}(c))})) \supseteq 
h^{-1}(h(a_{\id(g^{-1}(c))}))$ holds. Because we can choose 
an arbitrary numbering for the vertices of $G_A$, for an arbitrary 
country $d$ such that $h^{-1}(d) \neq \emptyset$ holds, we can 
assume that $a_1 \in h^{-1}(d)$. Combining this observation with 
the previous fact, we have that, for all $d \in V(G_M)$, if $h^{-1}(d) 
\neq \emptyset$ then there exists a country $c \in V(G_M)$ 
such that $h^{-1}(d) \subset g^{-1}(c)$. By the symmetry of the 
roles of $g$ and $h$, we have that $\{g^{-1}(c) | c \in V(G_M), 
g^{-1}(c) \neq \emptyset\}=\{h^{-1}(d) | d \in V(G_M), h^{-1}(d) 
\neq \emptyset\}$, and hence $h^{\prime}_l=h$, which means 
that $\phi_g \sim \phi_h$.
\end{Proofof}

Theorem\tume\ref{only_MERGER} has a corollary, which plays  
a key role with Theorem \ref{minimum_contraction} in our next 
algorithm for deciding t-equivalence. 
\begin{cor}\label{CONTRACT_ANYTIME}
Let $(G_A, G_M)$ be an arbitrary ordered pair of graphs, 
let $f$ be an arbitrary arrangement of $(G_A, G_M)$, let $g$ 
be an arbitrary $f$-irreducible arrangement, and let $h$ be 
an arbitrary arrangement which is t-equivalent to $f$. Then 
$g$ can be achieved from $h$ by a sequence of only SUBNET 
MERGERs.  
\end{cor}
\begin{Proof}
It is derived from Theorem\tume\ref{only_MERGER} and the fact that 
every $f$-irreducible arrangement is also $h$-irreducible. 
\end{Proof}

Theorem\tume\ref{minimum_contraction} tells us that 
two arrangements $f$ and $g$ are t-equivalent if and only if 
at least one $f$-irreducible configuration is equivalent to at 
least one $g$-irreducible configuration. 
Thanks to Corollary\tume\ref{CONTRACT_ANYTIME}, 
in order to find an $f$-irreducible arrangement, starting from 
the initial arrangement $f$, we can take an arbitrary proper 
SUBNET MERGER for the current arrangement and `contract'
the pebbles of its corresponding configuration, iteratively. 
By using Theorem\tume\ref{contact_decision}, we can find 
a proper SUBNET MERGER for a given arrangement in  
polynomial-time. If there exists no proper SUBNET MERGER, 
Corollary\tume\ref{CONTRACT_ANYTIME} tells us that the current 
arrangement is $f$-irreducible. If we obtain an $f$-irreducible 
arrangement $f_{\Irr}$ and an $g$-irreducible arrangement 
$g_{\Irr}$, by using Theorems\tume\ref{equivalence_decision_k_ge_2}
and \ref{equivalence_decision_k_eq_1}, we can check whether 
the corresponding two configurations $\phi_{f_{\Irr}}$ and 
$\phi_{g_{\Irr}}$ are equivalent or not in polynomial-time. 
Combining these observations, we have the following algorithm which 
decides in polynomial-time whether given two arrangements 
of an ordered pair of graphs are t-equivalent or not. 


\begin{alg}[t-equivalenceDecision($(G_A,G_M),f,g$)]\label{t-equivalence}
{\em ~\hfill
 \begin{namelist}{INSTANCEEE}
  \item[{\bf INPUT}]{Two arrangements $f,g$ of an ordered graph pair 
  $(G_A,G_M)$.} 
  \item[{\bf OUTPUT}]{Decide whether $f$ and $g$ are t-equivalent or not.}
 \end{namelist}
  \begin{alg-enumerate}
\item{Let $\phi_{f_{\Irr}}:=${\bf{IrreducibleConfiguration($(G_A, G_M),f$)}}
and \\
let $\phi_{g_{\Irr}}:=${\bf{IrreducibleConfiguration($(G_A,
    G_M),g$)}}.}\label{Make_Irreducible}
  \item{{\bf{If}} $\phi_{f_{\Irr}} \sim \phi_{g_{\Irr}}$ {\bf{then}} {\bf{return}} YES, 
otherwise {\bf{return}} NO.}\label{t-equivalence_decision}
  \end{alg-enumerate}}
\end{alg}

\vspace*{12pt}

\begin{alg}[IrreducibleConfiguration($(G_A,G_M),f$)]\label{IrrCon}
{\em ~\hfill
 \begin{namelist}{INSTANCEEE}
  \item[{\bf INPUT}]{An arrangement $f$ of an ordered graph pair 
  $(G_A,G_M)$.}
  \item[{\bf OUTPUT}]{An $f$-irreducible configuration $\phi_{f_{\Irr}}$.}
 \end{namelist}
  \begin{alg-enumerate}
  \item{Set $f_0:=f$. Make the isolated-agent set $\IA_0$ of $f_0$. 
Make the configuration $\phi_{0}$ associated with $f_0$.
Let $G_{0}$ denote the board graph of $\phi_{0}$, let $P_{0}$ 
denote the pebble set for $\phi_{0}$, and let $l_0$ denote 
the number of connected components of $G_{0}$.}\label{MakeFirstConfiguration} 
\item{{\bf{If}}  {\bf{ContractiblePair($\phi_{i}$)}}=$\{p,q\}$
{\bf{then}}  {\bf{goto}} \maru{\ref{next_arrangement}}; 
{\bf{else}} {\bf{goto}} \maru{\ref{end}}.}\label{Call_Oracle}  
\item{Set a new arrangement $f_{i+1}$ as 
\begin{equation*}
f_{i+1}(x):= \begin{cases}
            \phi_{i}(p), & \text{if $x \in p \cup q$}; \\
            f_{i}(x), & \text{if $x \notin p \cup q$}.
            \end{cases}
\end{equation*}}\label{next_arrangement}
\item{Make the isolated-agent set $\IA_{i+1}$ of $f_{i+1}$ 
by modifying the set $\IA_{i}$. Make the configuration 
$\phi_{i+1}$ associated with $f_{i+1}$ by modifying 
the configuration $\phi_{i}$.
Let $G_{i+1}$ denote the board graph of $\phi_{i+1}$, let $P_{i+1}$ 
denote the pebble set for $\phi_{i+1}$, and let $l_{i+1}$ denote 
the number of connected components of $G_{i+1}$.}\label{MakeConfiguration} 
\item{Set $i:=i+1$ and {\bf{goto}} $\maru{\ref{Call_Oracle}}$.}\label{increment}
\item{{\bf{Return}} the configuration $\phi_{i}$ 
associated with the arrangement $f_i$.}\label{end}
  \end{alg-enumerate}}
\end{alg}

\vspace*{12pt}

\begin{alg}[ContractiblePair($\phi_{i}$)]\label{Contractible_Pair}
{\em ~\hfill
 \begin{namelist}{INSTANCEEE}
  \item[{\bf INPUT}]{The isolated-agent set $\IA_i$  and the
  configuration $\phi_{i}$ defined in {\bf{Algorithm\tume\ref{IrrCon}}}.} 
  \item[{\bf OUTPUT}]{A subset $\{p,q\}$ of $P_{i} \cup \IA_i$ 
such that $p$ can contact $q$ and the graph $G_A[p \cup q]$ is 
connected, if any. Otherwise $\emptyset$.}
 \end{namelist}
\begin{alg-enumerate}
\item{{\bf{If}} there exists a pair of isolated-agents $\{a,b\} (\subseteq \IA_i)$ 
such that the pair $\{a,b\}$ is an edge of $G_A$ and that the pair 
$\{f_i(a), f_i(b)\}$ is an edge of  $G_M$, {\bf{then}} {\bf{return}}
$\{\{a\},\{b\}\}$.}
\item{Let $G_{i,j}(j=1,\ldots,l_i)$ denote each connected component 
of $G_{i}$. Let $P_{i,j}:={\phi_{i}}^{-1}(V(G_{i,j}))$ and let
$\phi_{i,j}$ denote the configuration of the pebble set 
  $P_{i,j}$ on the board graph $G_{i,j}$ as the restriction of
  $\phi_{i}$. 

{\bf{For}} $j:=1$ to $ l_i$ {\bf{do:}} 
 \begin{alg-enumerate}
  \item{Set $k_j:=|V(G_{i,j})| - |P_{i,j}|$, the vacancy size of 
  the configuration $\phi_{i,j}$. 
  Make the $k_j$-isthmus tree $T_{k_j}(G_{i,j})$ of the board graph 
  $G_{i,j}$ of the configuration $\phi_{i,j}$. Set $T:=T_{k_j}(G_{i,j})$.}
  \item{Choose a leaf $u$ of $T$. Let $B_u$ be a $k_j$-block of $G_{i,j}$
  corresponding to $u$. {\bf{If}} $T$ has a vertex $v$ such that 
  the length of the (unique) path of $T$ from $u$ to $v$ is $2$ 
  {\bf{then}} let ${_{u}I_{v}}$ denote the middle vertex of the path, 
  and let $B_v$ denote a $k_j$-block of $G_{i,j}$ corresponding to $v$ 
  {\bf{else}} $B_v:=\emptyset$. Let $\IC_i$ denote 
  the isolated-country set of $f_i$. And let $\IC_i[u]$ denote the set of
  all elements of $\IC_i$ adjacent (as vertices of $G_M$) to 
  at least one country of $B_u$. Let $P_{i}(u):=\phi_{i}^{-1}(B_u)$. 
  {\bf{If}} $G_M[B_u]$ is not a cycle graph, {\bf{then}} 
  $P_{i}(N_u):=\phi_{i}^{-1}(B_u \cup B_v \cup \IC_i[u])$ 
  {\bf{else}} $P_{i}(N_u):=\phi_{i}^{-1}(B_v \cup \IC_i[u])$ .
  {\bf{If}} there exists a pair $(p,q)$ of pebbles of $\phi_i$ such that 
  $p \neq q , p \in P_i(u), q \in P_i(N_u)$ and the agent subnetwork 
  $G_A[p \cup q]$ is connected, {\bf{then return}} $\{p,q\}$.}
  \item{{\bf{If}} $G_M[B_u]$ is a cycle graph, {\bf{and if}} 
the set $\phi_{i}^{-1}(B_u)$ contains a pair of distinct pebbles  
$\{p,q\}$ such that the agent subnetwork $G_A[p \cup q]$ is 
connected and that the pebble $p$ can contact the pebble $q$ 
along the cycle $G_M[B_u]$, {\bf{then}} {\bf{return}} $\{p,q\}$.}
  \item{{\bf{If}} $T$ has at least one edge {\bf{then}}
  {\bf{set}} $T:=T\setminus \{u, {_{u}I_{v}}\}$ and {\bf{goto}} 
  (\maru{\ref{find_pair}}-b).}
  \end{alg-enumerate}
}\label{find_pair}
\item{{\bf{Return}} $\emptyset$.}
\end{alg-enumerate}}
\end{alg}


\begin{thm}\label{t-equivalenceDecisionTheorem}
Algorithm\tume\ref{t-equivalence} works correctly. Its running time is 
$\Od(|E(G_M)|+(|V(G_M)|+|E(G_A)|)|V(G_A)|)$.
\end{thm}

\begin{Proof}
We have already explained (just before the description of the algorithm) 
the correctness of the algorithm. 

Now let us estimate the time complexity of the algorithm.  
Without loss of generality, here we assume that both the graphs 
$G_A$ and $G_M$ are connected. 
We can achieve the step \maru{\ref{t-equivalence_decision}} 
of Algorithm\tume\ref{t-equivalence} in $\Od(|E(G_A)|+|E(G_M)|)$-time 
by using Theorems\tume\ref{equivalence_decision_k_ge_2} and
\ref{equivalence_decision_k_eq_1}. 
The main body of our algorithm are Algorithm\tume\ref{IrrCon} 
and Algorithm\tume\ref{Contractible_Pair} (i.e. the step 
\maru{\ref{Make_Irreducible}} of Algorithm\tume\ref{t-equivalence}). 
The step \maru{\ref{MakeFirstConfiguration}} of
Algorithm\tume\ref{IrrCon} can be done in $\Od(|E(G_A)|+|E(G_M)|)$-time. 
We keep track of the family 
$S(G_i)$ of all the maximal isthmuses of the current board graph $G_i$ 
throughout Algorithm\tume\ref{IrrCon}, which is referred to each 
time when the $k$-isthmus trees of $G_i$ are made at the step 
\maru{\ref{find_pair}} in Algorithm\tume\ref{Contractible_Pair}. 
By using the list $S(G_i)$, we can keep down the total time of making 
all the $k_j$-isthmus trees $T_{k}(G_{i,j})$ $(j=1, \ldots, l_i)$ at the
step \maru{\ref{find_pair}} in Algorithm\tume\ref{Contractible_Pair} 
in $\Od(|V(G_M)|)$-time. The total running time of 
(\maru{\ref{find_pair}}-b) and (\maru{\ref{find_pair}}-c) 
in the step \maru{\ref{find_pair}} of
Algorithm\tume\ref{Contractible_Pair} is $\Od(|V(G_M)|+|E(G_A)|)$. 
Hence Algorithm\tume\ref{Contractible_Pair} can be done
in $\Od(|V(G_M)|+|E(G_A)|)$-time. Because the size of pebbles 
$|P_i|$ decreases after the step \maru{\ref{increment}} of 
Algorithm\tume\ref{IrrCon}, 
the steps \maru{\ref{Call_Oracle}}$-$\maru{\ref{increment}}
are executed at most $|V(G_A)|-1$ times. Hence the total running 
time of the {\bf{ContractiblePair}} oracle 
(i.e. Algorithm\tume\ref{Contractible_Pair}) during 
Algorithm\tume\ref{IrrCon} is $\Od((|V(G_M)|+|E(G_A)|)|V(G_A)|)$. 
The steps \maru{\ref{next_arrangement}} and
\maru{\ref{MakeConfiguration}} of Algorithm\tume\ref{IrrCon} 
can be done in $\Od(|V(G_A)|)$-time. Updating the list $S(G_i)$ 
after each step \maru{\ref{increment}} of Algorithm\tume\ref{IrrCon} 
also can be done in $\Od(|V(G_M)|)$-time. Hence the total running time 
of the steps \maru{\ref{Call_Oracle}}$-$\maru{\ref{increment}} with 
maintenance of the list $S(G_i)$ during Algorithm\tume\ref{IrrCon} is 
also $\Od((|V(G_M)|+|E(G_A)|)|V(G_A)|)$. Summing all the estimates in 
the above, the stated overall running time follows.
\end{Proof}

\begin{thm}\label{t-equivalenceDecision--The_Number_of_Transfers}
For an ordered pair $(G_A, G_M)$ of (general) graphs, and for two 
arrangements $f, g$ on the pair $(G_A, G_M)$, we can construct 
an explicit sequence of transfers from $f$ to $g$ of 
$\Eod(|V(G_M)|^2 \cdot |V(G_A)|)$-length. 
\end{thm}

\begin{Proof}
The total number of the {\bf{ContractiblePair}} oracle calles during 
Algorithm\tume\ref{IrrCon}, which is the same as the total number of 
contacts of `pebbles', is clearly at most $|V(G_A)|-1$. 
Kornhauser et al.\tume\cite{KMS1984} proved that the transitivity for 
the classical pebble motion problem can be done in $\Od(|V(G_M)|^2)$ 
moves and such a sequence of moves can be efficiently generated.  
By using these facts, we can construct an explicit sequence of transfers 
from $f$ to $g$ of $\Od(|V(G_M)|^2|V(G_A)|)$-length. 
On the other hand, Kornhauser et al.\tume\cite{KMS1984} also proved 
that the optimal transformation for the classical pebble motion problem 
requires $\Eod(|V(G_M)|^2 \cdot |V(G_A)|)$ moves, which guarantees that our 
optimal sequence of transfers also requires $\Eod(|V(G_M)|^2 \cdot |V(G_A)|)$-length. 
\end{Proof}


\section{Decision of Almightiness}\label{almighty}

In this section, we prove that the decision of almightiness
for AAP (and also for SGA)  is co-$\np$-complete. 

First we recall the definition of almightiness for AAP. 
Let $(G_A,G_M)$ be an ordered pair of simple undirected 
graphs. The pair $(G_A,G_M)$ is called almighty if the all 
arrangements on the pair $(G_A,G_M)$ are t-equivalent each 
other. 

Here let us assume that the input agent network $G_A$ is 
a connected graph, because this restriction may be a natural 
demand from several applications. Unfortunately, we will see 
that this restriction does not affect the time-complexity of 
our problem. 

\begin{prob}[The decision of almightiness]\label{Almightiness}
{\em ~\hfill
 \begin{namelist}{INSTANCEEEE}
  \item[{\bf INSTANCE}]{An ordered pair of simple undirected connected 
  graphs $(G_A,G_M)$.} 
  \item[{\bf PROBLEM}]{Is the ordered pair $(G_A,G_M)$ almighty?} 
 \end{namelist}}
\end{prob}

\begin{thm}\label{coNP}
The decision problem of almightiness for AAP is co-$\np$-complete.
\end{thm}

For example, suppose that $C$ is a circuit graph, $G_M$ is a graph 
whose complement $\ol{G_M}$ has a Hamiltonian circuit, and 
$|V(C)|=|V(G_M)|$ holds. Then there is an arrangement 
$f: V(C) \rightarrow V(G_M)$ of the ordered graph pair $(C, G_M)$ 
such that the $f$ is bijective and that every edge $\{u,v\} (\in E(C))$ is 
transformed into a non-edge $\{f(u),f(v)\} (\notin E(G_M))$. And hence 
in this case, the pair $(C, G_M)$ is not almighty. The Hamiltonian 
circuit problem (HC, for short) is a well known $\np$-complete 
problem. However we must not run away with idea that this example 
directly shows the co-$\np$-completeness of our problem, simply 
because the almightiness is so restrictive condition that, even if the 
complement of $G_M$ does not have a Hamiltonian circuit, there 
exist many types of candidate of such $G_M$ that the pair 
$(C, G_M)$ is not almighty. In stead of analyzing and classifying the 
complicated variety of the shapes of $G_M$ for which the pair 
$(C, G_M)$ is not almighty, let us take a strategy to provide
appropriate gadgets for and to attach them to each of the circuit 
graph $C$ and the input graph $G_M$. 

Before prove Theorem\tume\ref{coNP}, let us briefly analyze the 
time-complexity of the following restricted version of HC.

\begin{prob}[A restricted version of HC]\label{restricted_HC}
{\em ~\hfill
 \begin{namelist}{INSTANCEEEE}
  \item[{\bf INSTANCE}]{A simple undirected connected graph $H$ 
  whose complement $\ol{H}$ is also connected. } 
  \item[{\bf PROBLEM}]{Is the graph $H$ has a Hamiltonian circuit?} 
 \end{namelist}}
\end{prob}

\begin{la}\label{HCP_CONN}
The Problem\tume\ref{restricted_HC} is $\np$-complete. 
\end{la}

\begin{Proofof}{Lemma\tume\ref{HCP_CONN}}
It is clear that Problem\tume\ref{restricted_HC} is in $\np$. 
Now we will show a polynomial-time reduction from HC to 
Problem\tume\ref{restricted_HC}. Without loss of 
generality, we assume that the input graph $G$ of HC 
has at least $3$ vertices.  Let $v$ be an arbitrary vertex of $G$, 
and let $\{a,b,c\}$ be three vertices outside of $V(G)$. Let us 
construct a new graph $H$ by combining $G$ and the three 
vertices $\{a,b,c\}$ as follows: 
\begin{align*}
&V(H):=V(G) \cup \{a,b,c\};\\ 
&E(H):=E(G)\cup \big\{\{v,a\},\{a,b\},\{b,c\}\big\} \cup 
\big\{\{c,x\} : \{x,v\} \in E(G) \big\}.
\end{align*}
Then it is easy to see that the graph $H$ has a Hamiltonian circuit 
if and only if the original graph $G$ has a Hamiltonian circuit. 
Furthermore, the complement $\ol{H}$ of the graph $H$ is 
connected. Clearly, this reduction can be done in $\Od(|V(G)|)$-time. 
\end{Proofof}

\begin{Proofof}{Theorem\tume\ref{coNP}}
If an ordered pair $(G_A,G_M)$ of graphs is not almighty, as evidence 
for this fact, we can bring forward two arrangements $f$ and $g$ of 
the pair $(G_A,G_M)$ such that $f \not\cong g$ holds. Moreover, 
as we have shown in the section\tume\ref{equivalency}, 
we can verify the correctness of this evidence in polynomial-time. 
Hence Problem\tume\ref{Almightiness} is in co-$\np$. 

\noindent
Next, in order to prove the co-$\np$ completeness of 
Problem\tume\ref{Almightiness}, we will reduce the 
complement of Problem\tume\ref{restricted_HC} to 
Problem\tume\ref{Almightiness} in polynomial-time. 
Let $H$ be an input graph of Problem\tume\ref{restricted_HC}, 
$n:=|V(H)|$ be its order. Let $X$ be a set of $n$ vertices 
outside of $V(H)$, and let $\{a,b\}$ be a pair of distinct 
vertices outside of $V(H) \cup X$. By using these, we 
define the following new graph $G_A$:
\begin{align*}
&V(G_A):=V(H) \cup X \cup \{a,b\};\\ 
&E(G_A):=E( \ol{H} ) \cup {X \choose 2} \cup 
                 \big\{\{a,u\} : u \in V(H) \cup X\big\} 
                 \cup \big\{\{b,h\} : h \in V(H) \big\}.
\end{align*}
Then, let $C$ be a cycle graph such that its vertex-set is 
$V(C):=\{c_1,\ldots,c_n\}$ and that its edge-set is 
$E(C):=\big\{\{c_1,c_2\}, \ldots, \{c_{n-1},c_n\}, \{c_n,c_1\}\big\}.$ 
Let $Y$ be a set of $n$ vertices outside of $V(C)$, 
and let $\{p,q\}$ be a pair of distinct vertices outside of 
$V(C) \cup Y$. By using these, we define the following 
new graph $G_M$:
\begin{align*}
&V(G_M):=V(C) \cup Y \cup \{p,q\};\\ 
&E(G_M):=E(C) \cup {Y \choose 2}
                   \cup \Big( Y \times \big ( V(C) \cup \{q\} \big ) \Big ) 
                   \cup \big\{\{p,q\}\big\}.
\end{align*}

   \begin{center}
\includegraphics[width=10cm, angle=270]{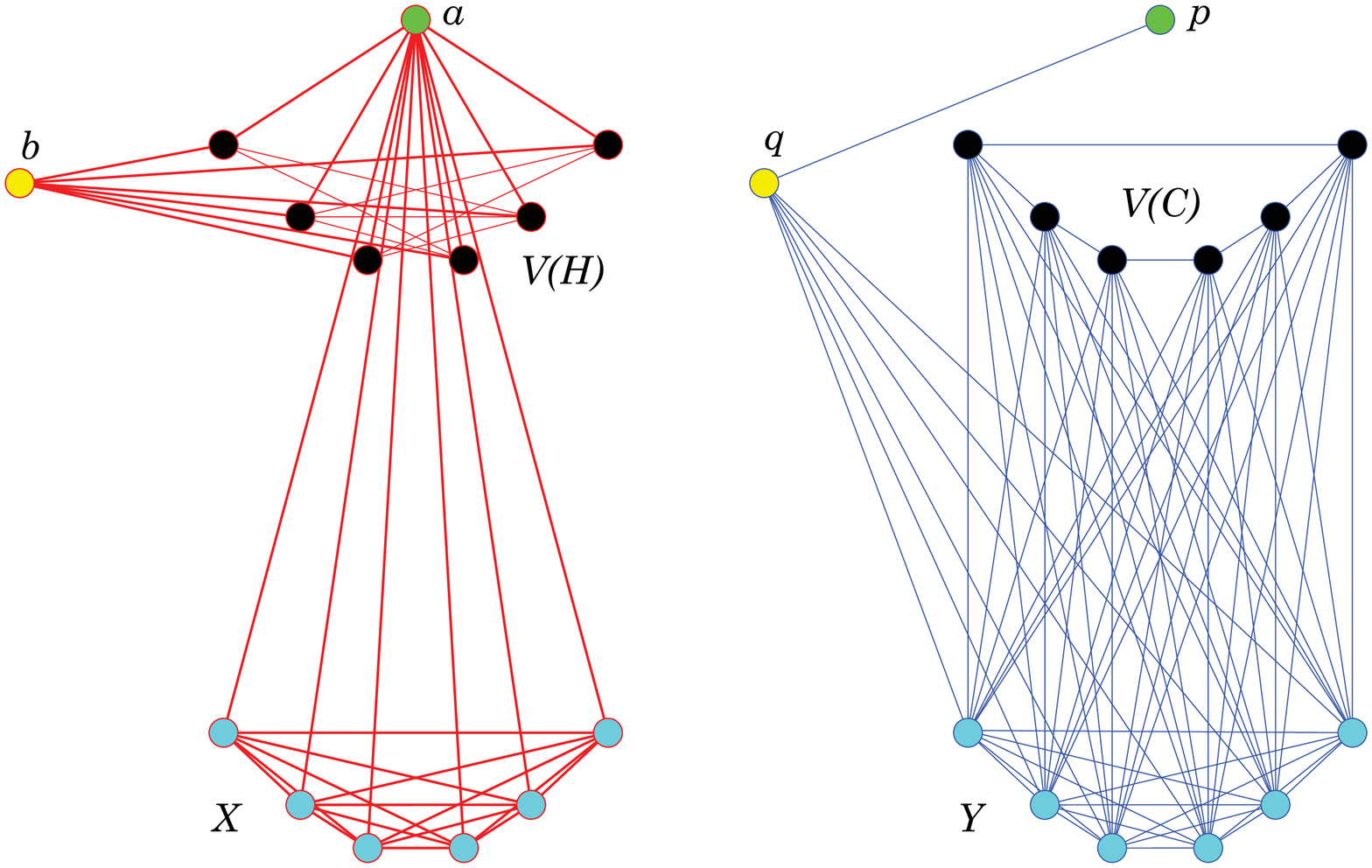}    
\textbf{Fig.3. }{The two graphs $G_A$ and $G_M$ in the proof of Theorem\tume\ref{coNP}.}
   \end{center}

Now suppose that $H$ has a Hamiltonian circuit. Then there exists 
a bijection $f: V(G_A) \rightarrow V(G_M)$ such that $f(a)=p$, $f(b)=q$, 
$f(X)=Y$, $f(V(H))=V(C)$ and, for every edge$\{u,u^{\prime}\}$ of 
the complement $\ol{H}$ of the graph $H$, 
the pair $\{f(u),f(u^{\prime})\}$ is a non-edge of $C$ (and hence of 
$G_M$).  This bijection $f$ can be thought as an arrangement of 
the pair $(G_A, G_M)$, and it is easy to see that this $f$ cannot be 
t-equivalent to any arrangement $g$ such that $\exists v \in V(H) 
\cup \{a,b\}, g(v) \neq  f(v)$. And hence the pair $(G_A, G_M)$ is 
not almighty.

On the contrary, let us assume that $H$ does not have any 
Hamiltonian circuit. We will show the fact that every arrangement $f$
of the pair $(G_A, G_M)$ is t-equivalent to the special 
arrangement $g: V(G_A) \to \{p\}$, which proves that 
the pair $(G_A, G_M)$ is almighty. 

First we prove the following technical lemma. 
\begin{la}\label{tec-lem}
For every vertex $v$ in $V(G_M) \setminus \{f(a)\}$, 
either $v \notin f(V(H) \cup \{b\})$ or $v \notin f(X)$ holds.
\end{la}
\begin{Proofof}{Lemma\tume\ref{tec-lem}}
The graph $G_A[f^{-1}(v)]$ must be connected and
every connected subgraph of $G_A -a$ is a subgraph
of either $G_A[\{b\} \cup V(H)]$ or $G_A[X]$, which proves the lemma.  
\end{Proofof}

\setcounter{case}{0}

We divide our proof into the several cases. 
\begin{case}{\em
$f(a) \in Y.$ 

\begin{subcase}{\em There exists a vertex $h$ in $V(H)$ such that 
$f(h) \in V(G_M) \setminus \{p\}$. 

\noindent
Because the pair $\{a,h\}$ is an edge of $G_A$ and 
the pair $\{f(a),f(h)\}$ is an edge of $G_M$, we may transfer 
the graph $G_A[f^{-1}(f(h))]$ onto the vertex $f(a)$ of $G_M$ and 
merge it into $G_A[f^{-1}(f(\{a,h\}))]$. Next, because 
every vertex of $G_A$ is adjacent to at least one of $\{a,h\}$, 
we can move all the vertices of $f^{-1}(V(G_M) \setminus \{p\})$ 
onto the vertex $f(a)$ of $G_M$. Then,
$G_A[f^{-1}(V(G_M) \setminus \{p\})]$ is on the vertex $f(a)$ of 
$G_M$ and all the other vertices of $G_A$ are on the vertex $p$ 
of $G_M$. Finally, we may transfer the graph 
$G_A[f^{-1}(V(G_M) \setminus \{p\})]$ onto the vertex $p$ and 
merge it into $G_A$ on $p$. Hence the arrangement $f$ of 
the pair $(G_A, G_M)$ is t-equivalent to the special arrangement 
$g: V(G_A) \to \{p\}$.}
\end{subcase}

\begin{subcase}{\em $f(V(H))=\{p\}, f(b)=q$. 

\noindent
Since every vertex  of $V(H)$ is adjacent to $b$ and the pair $\{p,q\}$
is an edge of $G_M$, we can transfer the graph $G_A[f^{-1}(p)]$ onto 
the vertex $q$ of $G_M$ and merge it into $G_A[f^{-1}(\{p,q\})]$, which 
means that the case is reduced to {\bf Subcase 1.1}. }
\end{subcase}

\begin{subcase}{\em $f(V(H))=\{p\}, f(b) \neq q$. 

\noindent
Since every vertex in $f^{-1}(q)$ is adjacent to the vertex $a$ 
and the pair $\{q,f(a)\}$ is an edge of $G_M$, we may transfer 
the graph $G_A[f^{-1}(q)]$ onto the vertex $f(a)$ of $G_M$ and 
merge it into $G_A[f^{-1}(\{q,f(a)\})]$. Then the vertex $q$ will be 
unoccupied, and hence we can transfer the graph $G_A[f^{-1}(p)]$ 
onto the unoccupied vertex $q$, which means that the case is 
reduced to {\bf Subcase 1.1}. } 
\end{subcase}
}
\end{case}

\begin{case}{\em
$f(a) \in V(C).$

\noindent
Note that Lemma\tume\ref{tec-lem} tells us that 
$p \notin f(V(H)) \cap f(X)$ holds. 

\begin{subcase}{\em
There exists a vertex $u$ in $V(H) \cup X$ such that $f(u) \in Y$.

\noindent
In this case, since the pair $\{a,u\}$ is an edge of $G_A$ and 
the pair $\{f(a),f(u)\}$ is an edge of $G_M$, we can transfer 
the graph $G_A[f^{-1}(f(a))]$ onto the vertex $f(u)$ of $G_M$ and 
merge it into $G_A[f^{-1}(f(\{a,u\}))]$, which means that the 
case is reduced to {\bf Case 1}.} 
\end{subcase}

\begin{subcase}{\em $f(V(H) \cup X) \cap Y=\emptyset$ and 
there exist a vertex $w$ in $V(C)\cup\{q\}$ such that $|f^{-1}(w)| \geqq 2$. 

\noindent
In this case, we can regard 
this connected subgraph $G_A[f^{-1}(w)]$ as if a single pebble 
on the board graph $G_M$. Because $b$ is a unique vertex of 
$G_A$ which is not adjacent to $a$, at least one vertex $y$ of $Y$ 
is unoccupied. We can move the pebble $G_A[f^{-1}(w)]$ onto this 
unoccupied vertex $y$, and hence the case is reduced to either 
{\bf Case 1} or {\bf Subcase 2.1}.}
\end{subcase}

\begin{subcase}{\em $f(V(H) \cup X) \cap Y = \emptyset$ and 
every vertex $w$ in $V(C)\cup\{q\}$ satisfies $|f^{-1}(w)| \leqq 1$.

\noindent
In this case, because $p \notin f(V(H) \cup \{b\}) \cap f(X)$, 
the Pigeonhole principle shows that the set $X \cup V(H)$ includes
a vertex $u$ such that $\{f(a), f(u)\}$ is an edge of the cycle graph
$C$. Since $\{a,u\}$ is an edge of $G_A$, we can merge $a$ and $u$ 
into the graph $G_A[\{a,u\}]$ and put it on the vertex $f(a)$ of $G_M$. 
Moreover, we can regard this $G_A[\{a, u\}]$ as if a single pebble 
on the board graph $G_M$. Since at least one vertex $y$ of $Y$ is 
unoccupied, we can move this pebble $G_A[\{a, u\}]$ onto 
the unoccupied vertex $y$ of $G_M$, which means that the case is 
reduced to {\bf Case 1}. }
\end{subcase}

}
\end{case}

\begin{case}{\em 
$f(a)=q.$

\begin{subcase}{\em The set $V(H) \cup X$ includes a vertex $u$ 
such that $f(u) \in Y$. 

\noindent
The pair $\{a,u\}$ is an edge of $G_A$ and the pair 
$\{f(a),f(u)\}$ is an edge of $G_M$. Hence we can transfer 
the graph $G_A[f^{-1}(f(a))]$ onto the vertex $f(u)$ of
$G_M$ and merge it into $G_A[f^{-1}(f(\{a,u\}))]$, 
which means that the case is reduced to {\bf Case 1}. }
\end{subcase}

\begin{subcase}{\em $f(V(H) \cup X) \cap Y = \emptyset$ and 
the set $V(H) \cup X$ includes a vertex $u$ such that $f(u)=p$. 

\noindent
In the same way as the previous case, the pair $\{a,u\}$ is an edge of 
$G_A$ and the pair $\{f(a),f(u)\}$ is an edge of $G_M$. 
Hence we can transfer the graph $G_A[f^{-1}(f(u))]$ onto the vertex 
$q$ of $G_M$ and merge it into $G_A[f^{-1}(\{p,q\})]$. 
Because $b$ is a unique vertex of $G_A$ which is not adjacent to $a$, 
at least one vertex $y$ of $Y$ is unoccupied. We can transfer the graph 
$G_A[f^{-1}(\{p,q\})]$ onto this unoccupied vertex $y$ of $G_M$, and 
hence the case is reduced to {\bf Case 1}. }
\end{subcase}

}
\end{case}

\begin{case}{\em 
$f(a)=p.$

\begin{subcase}{\em $f(b)=p$. 

\noindent
Since every vertex of $G_A$ is adjacent to at least one of the pair 
$\{a,b\}$, we can transfer $G_A[f^{-1}(p)]$ onto the vertex $q$ of $G_M$ 
and merge it into $G_A[f^{-1}(\{p,q\})]$, which means that the case is 
reduced to {\bf Case 3}. }
\end{subcase}

\begin{subcase}{\em There exists a vertex $u$ in $V(H) \cup X$ such that 
$f(u)=q$. 

\noindent
The pair $\{a,u\}$ is an edge of $G_A$ and the pair $\{f(a),f(u)\}$ is
an edge of $G_M$. Hence we can transfer the graph $G_A[f^{-1}(f(a))]$ 
onto the vertex $q$ of $G_M$ and merge it into $G_A[f^{-1}(\{p,q\})]$,  
which means that the case is reduced to {\bf Case 3}. }
\end{subcase}

\begin{subcase}{\em $q \notin f(V(H) \cup X)$ and 
$f(b) \neq p$ and there exists a vertex $h$ of $V(H)$ such that 
$f(h) \in Y$. 

\noindent
Because $f(h) (\in Y)$ is adjacent to every vertex of $V(C) \cup Y \cup
\{q\}$ and $\{b,h\}$ is an edge of $G_A$, we can transfer the graph 
$G_A[f^{-1}(f(b))]$ onto the vertex $f(h)$ of $G_M$ and merge it into 
$G_A[f^{-1}(f(\{b,h\}))]$. Moreover we can also transfer the
graph $G_A[f^{-1}(f(\{b,h\}))]$ onto the vertex $q$ of $G_M$ 
because either $f^{-1}(q)=\emptyset$ or $f(b)=q$ holds. Hence 
the case is reduced to {\bf Subcase 4.2}. }
\end{subcase}

\begin{subcase}{\em $q \notin f(V(H) \cup X)$ and 
$f(b) \in Y$ and $f(V(H)) \subseteq V(C)$. 

\noindent
For every vertex $h$ in $V(H)$, the pair $\{h,b\}$
is an edge of $G_A$ and the pair $\{f(h),f(b)\}$ is an edge of $G_M$, 
and hence we can transfer the graph $G_A[f^{-1}(f(h))]$ 
onto the vertex $f(b)$ of $G_M$ and merge it into $G_A[f^{-1}(f(\{b,h\}))]$,  
which means that the case is reduced to {\bf Subcase 4.3}.}
\end{subcase}

\begin{subcase}{\em $q \notin f(V(H) \cup X)$ and 
$f(b) \notin Y$ and $f(V(H)) \subseteq V(C)$. 

\noindent
In this case, we have that $f^{-1}(Y) \subseteq X$, 
and hence, for every vertex $x$ in $f^{-1}(Y)$, we may 
transfer the graph $G_A[f^{-1}(Y)]$ onto the vertex $f(x)$ of
$G_M$. Because $|Y| = n \geqq 2$, $Y$ includes at least one vertex $y$ 
other than $f(x)$. We may assume that this vertex $y$ will be unoccupied. 
Because $f(V(H)) \subseteq V(C)$, and because $H$ does not have any 
Hamiltonian circuit, $V(H)$ includes a pair of distinct vertices
$\{h_1,h_2\}$ such that either the pair $\{f(h_1), f(h_2)\}$ is an edge 
of $C$ or $f(h_1)=f(h_2)$ holds. Hence we may transfer the graph 
$G_A[f^{-1}(f(h_1))]$ onto the vertex $f(h_1)$ of $G_M$ and merge 
it into $G_A[f^{-1}(f(\{h_1, h_2\}))]$. Next, we may transfer 
the graph $G_A[f^{-1}(f(\{h_1, h_2\}))]$ onto the unoccupied vertex 
$y$ of $G_M$, which means that the case is reduced to {\bf Subcase 4.3}. }
\end{subcase}
}
\end{case}
Now our proof is completed. 
\end{Proofof}

\begin{cor}\label{SGA_NP-complete}
The decision problem of almightiness for SGA is co-$\np$-complete.
\end{cor}
Its proof is almost the same as the one of Theorem\tume\ref{coNP} 
and will be omitted. (Please use the one vertex $v_X$ (resp. $v_Y$) 
except for the gadget $X$ (resp. $Y$).) 


%

%
\end{document}